\documentclass[12pt]{article}

\title{Cone distribution functions and quantiles for multivariate random variables}

\author{
Andreas H. Hamel\footnote{Free University Bozen, Faculty of Economics and Management, \href{mailto:andreas.hamel@unibz.it}{andreas.hamel@unibz.it}}, Daniel Kostner\footnote{Free University Bozen, Faculty of Economics and Management, \href{mailto:daniel.kostner@economics.unibz.it}{daniel.kostner@economics.unibz.it}}
}
\date{{\small \today}}

\usepackage{array} 
\usepackage{verbatim} 
\usepackage{float, graphicx}

\usepackage{amsmath, amssymb, enumerate, color}

\newtheorem{theorem}{Theorem}
\newtheorem{corollary}[theorem]{Corollary}
\newtheorem{remark}[theorem]{Remark}
\newtheorem{lemma}[theorem]{Lemma}
\newtheorem{definition}[theorem]{Definition}
\newtheorem{proposition}[theorem]{Proposition}
\newtheorem{example}[theorem]{Example}

\numberwithin{equation}{section}  
\numberwithin{figure}{section}    
\numberwithin{table}{section}     
\numberwithin{theorem}{section}

\newcommand{\of}[1]{\ensuremath{\left( #1 \right)}}

\newcommand{\cb}[1]{\ensuremath{ \left\{ #1 \right\} }}
\newcommand{\sqb}[1]{\ensuremath{ \left[ #1 \right] }}

\newcommand{\bs}{\backslash}

\newcommand{\pend}{ \hfill $\square$ \medskip}

\newcommand{\A}{\ensuremath{\mathcal{A}}}

\newcommand{\G}{\ensuremath{\mathcal{G}}}

\newcommand{\R}{\mathrm{I\negthinspace R}}

\newcommand{\N}{\mathrm{I\negthinspace N}}

\newcommand{\cl}{{\rm cl \,}}

\newcommand{\bd}{{\rm bd \,}}
\newcommand{\co}{{\rm co \,}}

\newcommand{\qint}{{\rm qint \,}}

\newcommand{\Int}{{\rm int\,}}

\newcommand{\One}{\mathrm{1\negthickspace I}}

\usepackage[
colorlinks=true, linkcolor=blue, citecolor=blue, urlcolor=blue, filecolor=blue
auskommentieren und unten pdfborder aktivieren]{hyperref}

\definecolor{color0}{gray}{.50}
\definecolor{color1}{rgb}{0,.2,.8}
\definecolor{color2}{rgb}{1,.2,0}
\definecolor{color3}{rgb}{.8,.5,1}

\begin{document}
\maketitle

\abstract{Set-valued quantiles for multivariate distributions with respect to a general convex cone are introduced which are based on a family of (univariate) distribution functions rather than on the joint distribution function. It is shown that these quantiles enjoy basically all the properties of univariate quantile functions. Relationships to families of univariate quantile functions and to depth functions are discussed. Finally, a corresponding Value at Risk for multivariate random variables as well as stochastic orders are introduced via the set-valued approach.}


\section{Introduction}

When it comes to quantiles for multivariate random variables, there is no ``silver bullet," but several very different approaches. Some are based on the joint distribution function e.g.  \cite{EmbrechtsPuccetti06JMA}, others on different statistical depth functions e.g. \cite{HallinEtAl10AS, KongMizera12StSi, ZuoSerfling00AS}. While the former involve the ``natural" ordering cone $\R^d_+$, the latter rarely involve any ordering for the values of the random variable. The joint distribution approach is sometimes coupled with copulas which leads to real-valued quantiles for multivariate distribution e.g. via Kendall's distribution as in \cite[Definition 5]{SalvadoriMicheleDurante11HESS}. 

The following question remains. What is the (upper and lower) quantile of a multivariate variable if the decision maker/analyst has a preference for the data points which is, for example, a vector order in the outcome space? This order relation does not enter the picture through statistical analysis, it is rather a given object that should influence statistical procedures: If the order changes, the corresponding quantiles, outlyingness notion, stochastic orders and risk measures should also change. Such an order occurs frequently and naturally for financial data in the presence of transaction costs. If the latter are proportional, the order is generated by the so-called solvency cone (which is usually different from $\R^d_+$), see e.g. \cite{HamelHeyde10SIFIN}, \cite{HamelHeydeRudloff11MFE} and the references therein.

An attempt to incorporate general orders into the statistical analysis of multivariate data is one topic of \cite{BelloniWinkler11AS} which led to concepts and formulas which are "very far" from the univariate case. Belloni and Winkler write (p. 1126) `The fundamental difficulty in reaching agreement on a suitable generalization of univariate quantiles is arguably the lack of a natural ordering in a multidimensional setting.' Therefore, the main goal of the present paper is to develop a theory which can deal with any (vector) order in the outcome space and still runs completely parallel to the univariate case.

The same remarks apply to the definition of the Value at Risk (VaR) for multivariate positions. VaR for univariate variables--just a quantile--is a common risk evaluation tool in finance and insurance which has its advantages and drawbacks. Among the former certainly is that VaR can be used to define new and sometimes more appropriate risk measures such as the Average (or Conditional) Value at Risk (AVaR). An example for a multivariate VaR can be found in \cite{EmbrechtsPuccetti06JMA} which is also based on the joint distribution function while \cite{CascosMolchanov07FS} gives a version which is related to general order relations (but not completely sound).

Undoubtedly, there is demand for quantile-like concepts in multivariate analysis \cite[p. 1125]{BelloniWinkler11AS}: `Naturally, the quantiles of a multivariate random variable are also of interest, and the search for a multidimensional counterpart of the quantiles of a random variable has attracted considerable attention in the statistical literature.' It has been remarked \cite[p. 214]{Serfling02SN} that `various ad hoc quantile-type multivariate methods have been formulated, some vector-valued in character, some univariate, and the term ``quantile" has acquired rather loose usage.' This sparked several axiomatic approaches to statistical depth functions and depth regions e.g. \cite{ZuoSerfling00AS, Serfling02SN} and also \cite{CascosMolchanov07FS} with financial risk measures in view.

The introduction of quantile-like concepts on the one hand and depth (regions) on the other hand often comes with `sharp methodological differences'  (\cite[p. 636]{HallinEtAl10AS}) described as follows in the same reference: `While quantiles resort to analytical characterizations through inverse distribution functions or $L_1$ optimization, depth often derives from more geometric considerations such as halfspaces, simplices, ellipsoids and projections.' The quantiles obtained through $L^1$ optimization techniques such as in \cite{Chaudhuri96JASA} often rely on an indexing procedure e.g. ``by elements of the open unit ball" (\cite[p. 863]{Chaudhuri96JASA}, see also \cite{HallinEtAl10AS}) in the outcome space and do not involve a (vector) order for the data points.

In this note, we propose a novel approach which could be seen as an attempt to bridge analytic and geometric concepts and arrives at formulas which are as close to the univariate case as possible. The approach rests on recent developments in set optimization and set-valued variational analysis as surveyed in \cite{HamelEtAl15Incoll} and admits to involve a general vector order for the values of the multivariate variable. The one major fact one has to cope with--different from the univariate case--is that quantiles become functions mapping into well-defined lattices of sets and thus are set-valued in nature. 

We start by introducing a generalization of the Tukey halfspace depth function which we call cone distribution function since, on the one hand, it behaves pretty much like a (joint) distribution function and, on the other hand, depends on a cone which can be (very) different from $\R^d_+$. Next, set-valued quantile functions for multivariate variables are introduced with the following features: (1) set-valued lower and upper quantiles are basically set-valued inverses of the cone distribution function and its ``strict" counterpart, (2) they produce functions with values in (two different) complete lattices of sets, (3) the ``set-valued" formulas can be understood completely analogous to the univariate case, (4) any vector order can be dealt with, (5) there is no need for an indexing procedure or the choice of a direction during the statistical analysis, our quantiles only depend on the given data and the order.

Finally, our new concepts are applied in order to define a set-valued Value at Risk and a stochastic dominance order for multivariate variables. 

Our constructions provide evidence that there are two ``natural ordering(s) in a multivariate setting" on a set level: one for lower and one for upper quantiles. This point of view is supported by the strong link to set optimization which is explained in Section \ref{SecSOP}. Moreover, our approach also resolves the ambiguity which stems from the fact that the (joint) cdf and the (joint) survival function produce different concepts (usually called "lower orthant" and "upper orthant" quantile, Value at Risk, stochastic order etc. as in \cite{EmbrechtsPuccetti06JMA} and \cite{MuellerStoyan02Book}, for example) due to the non-totalness of vector orders in higher dimensional spaces.

Examples are given, and comparisons with existing concepts from the literature conclude the paper along with a discussion of desirable generalizations and extensions. For the convenience of the reader, an appendix with basic concepts related to convex cones is added.

\section{Distribution functions associated to a cone}

Let $\of{\Omega, \mathcal F, P}$ be a probability space and $X \colon \Omega \to \R^d$ a multivariate random variable, i.e. an $(\mathcal F, \mathcal B^d)$-measurable function. A standard concept in probability theory and statistics is the joint distribution function $F^{jdf}_X \colon \R^d \to [0,1]$ defined by
\[
F^{jdf}_X(z) = P\of{X \in z - \R^d_+}.
\]
It involves the component-wise (partial) order in $\R^d$ generated by the closed convex cone $\R^d_+$. A natural idea would be to replace $\R^d_+$ by a general convex cone $C \subseteq \R^d$ with $0 \in C$ and define an analog to $F^{jdf}_X$. We will not follow this path for several reasons. One of them is that the resulting quantile (set) is non-convex in general (see Example \ref{Joint} below).

Instead, we propose to base the discussion upon different objects which might replace the joint distribution function for some purposes, namely a family of (ordinary) cumulative distribution functions and a cone distribution function which turns out to be different from the joint distribution function even if the cone is $\R^d_+$. Similar functions have already been considered in \cite{RousseeuwRuts99Met, HallinEtAl10AS, KongMizera12StSi}.

We recall a few concepts from the theory of ordered vector spaces in order to fix the notation. For $w \in \R^d\bs\{0\}$, the set $H^+(w) = \cb{z \in \R^d \mid w^Tz \geq 0}$ is the closed homogeneous halfspace with normal $w$. A set $C \subseteq \R^d$ is called a cone if $s>0$, $z \in C$ imply $sz \in C$. A cone $C$ is a convex set if, and only if, it is closed under addition, i.e. $x,y \in C$ implies $x + y \in C$, and in this case it is called a convex cone. A convex cone $C \subseteq \R^d$ with $0 \in C$ generates a vector preorder $\leq_C$ by means of
\[
z \leq _C z' \quad \Leftrightarrow \quad z' - z \in C.
\]
This means that $\leq_C$ is a reflexive and transitive relation which is compatible with the algebraic operations of the linear space $\R^d$. Vice versa, every such preorder $\leq$ can be represented by the convex cone $C_\leq = \cb{z \in \R^d \mid 0 \leq z}$. The preorder $\leq_C$ is antisymmetric, i.e. a partial order, if and only if, $C \cap -C = \cb{0}$. In this case, $C$ is called pointed, but we will not assume in the following that $C$ is pointed. On the contrary, the case $C = H^+(w)$ for some $w \in \R^d\bs\{0\}$ is a valid option.

If $C \subseteq \R^d$ is a convex cone, the set $C^+ = \cb{w \in \R^d \mid \forall z \in C \colon w^Tz \geq 0}$ is called its (positive) dual cone (sometime also polar cone). The dual of a convex cone always is a closed convex cone. The dual cone $C^+$ is said to have a base $B^+$ if $B^+ \subseteq C^+$ is closed convex set which does not contain $0 \in \R^d$ such that for each $w \in C^+\bs\{0\}$ there exist unique $s>0$, $b \in B^+$ with $w = sb$. 

The bipolar theorem gives a relationship between a cone $C$ and its dual $C^+$. It states that a convex cone $C$ is closed if, and only if, 
\[
C = \cb{z \in \R^d \mid \forall w \in C^+ \colon w^Tz \geq 0}.
\]
This means, $C$ is a closed convex cone if, and only if, $C = C^{++} := (C^+)^+$. In this case, the relation $\leq_C$ has a representation by a family of scalar functions, i.e.
\[
z \leq_C y \quad \Leftrightarrow \quad \forall w \in C^+ \colon w^Tz \leq w^Ty \quad \Leftrightarrow \quad \forall w \in C^+ \colon z \leq_{H^+(w)}y.
\]
This means that $\leq_C$ is represented as intersection of total orders generated by the closed halfspaces $H^+(w)$ for $w \in C^+\bs\{0\}$.

If $C^+$ has a base $B^+$ then, of course, it is enough to let $w$ run through $B^+$ instead of $C^+$ in the above representation of $\leq_C$.

If $C = \{0\}$, then $C^+ = \R^d$, and there is no base for $C^+$. If $C = \R^d_+$, then $C^+ = \R^d_+$, and the set $B^+ = \cb{w \in \R^d_+ \mid w^Te = 1}$ is a base of $\R^d_+$ where $e = (1, 1, \ldots, 1)^T \in \R^d$. If $C = H^+(w)$ for some $w \in \R^d\bs\{0\}$, then $C^+ = \cb{sw \mid s \geq 0}$, and $B^+ = \{w\}$ is a base of $C^+$. In particular, if $d = 1$, $C = \R_+$, then $B^+ = \{1\}$ is a base of $C^+ = \R_+$, and this simple device produces the scalar special case in all the considerations in this note.

The following definition is the departing point for defining set-valued quantiles. 

\begin{definition}
\label{DefDistFuncFam}
The set $\cb{F_{X,w}}_{w \in C^+\bs\{0\}}$ of functions $F_{X,w} \colon \R^d \to [0,1]$ defined by
\[
F_{X,w}(z) = F_{w^TX}(w^Tz) = P\of{X \in z - H^+(w)}
\]
is called the family of (cumulative) distribution functions for $X \colon \Omega \to \R^d$ with respect to $C$. The function $F_{X, C} \colon \R^d \to [0,1]$ defined by
\[
F_{X, C}(z) = \inf_{w \in C^+\bs\{0\}} F_{X,w}(z) 
\]
is called cone distribution function of $X$ (with respect to the cone $C$) or just $C$-distribution function.
\end{definition}

Note that $F_{X, 0}(z) \equiv 1$, so $w = 0$ can be excluded in the definition of $F_{X, C}(z)$. Clearly, $F_{X,sw}(z) = F_{X,w}(z)$ for all $s>0$, so if $C^+$ has a base $B^+$, then it even suffices to consider $\cb{F_{X,w}}_{w \in B^+}$, and in this case $F_{X, C}(z) = \inf_{w \in B^+} F_{X,w}(z)$. 

Assume that for $z \in \R^d$ the infimum defining the $C$-distribution function is attained at $\bar w \in C^+\bs\{0\}$. Then $F_{X, C}(z) = P\of{X \in z - H^+(\bar w)}$, which means that the halfspace $z -  H^+(\bar w)$ with $z$ at its boundary is least likely to contain values of $X$ among all halfspaces with normals in $C^+\bs\{0\}$. There is, of course, a strong link to Tukey depth functions as explained in Section \ref{SecConnectTukey} below. For an empirical version, compare e.g. \cite[formula (1.4)]{StruyfRousseeuv99JMA}.

Finally, since $C^+ = \of{\cl C}^+$ the $C$-distribution function coincides with the corresponding $\cl C$-distribution function. Therefore, there is no loss in generality by assuming that $\emptyset \neq C \neq \R^d$ is a closed convex cone. This is a standing assumption from now on.

\begin{remark}
\label{RemCDFvsJDF}
It is easy to see that $F_{X, C}(z) \geq P(X \in z - C)$ for all $ z \in \R^d$. This inequality is strict in general even if $C = \R^d_+$ (the case of the joint distribution function) as one may already observe for the bivariate standard normal distribution: In this case, $F^{jdf}_X(0) = P(X \in -\R^d_+) = \frac{1}{4}$ while $P(X \in -H^+(w)) = \frac{1}{2}$ for all $w \in \R^d_+\bs\{0\}$, so $F_{X, C}(0) = \frac{1}{2}$. One may also observe that this phenomenon is related to the non-totalness of the order $\leq_C$, i.e. the existence of non-comparable elements. On the other hand, if $C = H^+(\bar w)$ for some $\bar w \in \R^d\bs\{0\}$, then $F_{X, C}(z) =  F_{X, \bar w}(z) = P(X \in z - C)$ for all $z \in \R^d$.
\end{remark}

\begin{example}
Let $w = e^i = (0, \ldots, 0,1,0, \ldots, 0)^T$ with i-th component equal to 1. Then
\[
\forall z \in \R^d \colon F_{X,w}(z) = F_{w^TX}(w^Tz) = P\of{X_i \leq z_i}
\]
which is the marginal distribution function of $X$ with respect to the i-th component for $i \in \cb{1, \dots, d}$. However, it might very well happen that $e^i \not\in C^+$ for a cone $C$ and some $i \in \cb{1, \ldots, d}$ in which case the corresponding marginal distribution does not seem to be a relevant object in our framework.
\end{example}

\begin{example}
\label{Ex1}
By a slight abuse of notation, let $X$ and $Y$ be two independent random variables, uniformly distributed on $(0,1)$. Then, the joint distribution function of the bivariate random variable $(X, Y)$ is
\begin{align*}
F^{jdf}_{X, Y}\of{x,y} = & \, P\of{(X,Y) \in (x,y) - \R^2_+} \\
& = \begin{cases}
    0, & \text{for} \quad y <0 \quad \text{or} \quad x <0, \\
	xy, & \text{for} \quad 0 \leq y \leq 1, \quad 0 \leq x \leq 1,\\
	y, & \text{for} \quad x > 1, \quad 0 \leq y \leq 1, \\
	x, & \text{for} \quad y > 1, \quad 0 \leq x \leq 1, \\
	1, & \text{for} \quad x > 1, \quad y > 1.
\end{cases}
\end{align*}
whereas the bivariate lower $\R^2_+$-distribution function of $(X,Y)$ is
\begin{align*}
F_{(X,Y), \R^2_+}\of{(x,y)} = & \inf_{w \in \R^2_+\bs\{0\}} P\of{(X,Y) \in (x,y) - H^+(w)}  \\
& = \begin{cases}
    0, & \text{for} \quad y <0 \quad \text{or} \quad x <0, \\
	 \min\cb{x, y}, & \text{for} \quad 0 \leq y \leq 1, \quad 0 \leq x \leq 1,\\
	y, & \text{for} \quad x > 1, \quad 0 \leq y \leq 1, \\
	x, & \text{for} \quad y > 1, \quad 0 \leq x \leq 1, \\
	1, & \text{for} \quad x > 1, \quad y > 1.
\end{cases}
\end{align*}
The expression for the case $(x,y) \in [0,1] \times [0,1]$ follows since $B^+ = \cb{w \in R^2 \mid w_1 + w_2 = 1}$ is a base for $\R^2_+$ and
\[
\inf_{w \in \R^2_+\bs\{0\}} \frac{w_1x + w_2y}{w_1 + w_2} = \inf_{w \in B^+} w_1x + w_2y = \min\cb{x, y}.
\]
Moreover, $xy \leq \min\cb{x, y}$ with strict inequality for $(x,y) \in (0,1) \times (0,1)$. Hence $F_{X, Y}\of{x,y} <  F_{(X,Y), \R^2_+}\of{(x,y)}$ on $(0,1) \times (0,1)$.
\end{example}

A few elementary properties of cone distribution functions are collected in the following result which needs one more concept related to cones: The set $\qint C = \cb{z \in C \mid \forall w \in C^+\bs\{0\} \colon w^Tz > 0}$ is called the quasi-interior of $C$. Moreover, if $\leq_C$ is applied to $\R^d$-valued random variables it is understood in an almost sure sense.

\begin{proposition}
\label{PropCDistribution} The cone distribution function $F_{X, C}$  has the following properties:

(a) It is {\em affine equivariant}, i.e. if $b \in \R^d$ and $A \in \R^{d \times d}$ is an invertible matrix, then
\[ 
\forall z \in \R^d \colon F_{AX + b, AC}\of{Az + b} = F_{X, C}\of{z}.
\]

(b) It is a {\em monotone non-decreasing} function of $z$ with respect to $\leq_C$, i.e. if $y \leq_{C} z$, then $F_{X, C}\of{y} \leq F_{X, C}\of{z}$. 

(c) It is a {\em monotone non-increasing} function of $X$ with respect to $\leq_C$, i.e. if $X \leq_{C} Y$, then $F_{X, C}\of{z} \geq F_{Y, C}\of{z}$ for all $z \in \R^d$.

(d) It is {\em right-continuous}, i.e. if $\{z_n\}_{n \in \N} \subseteq \R^d$  a sequence with $\lim_{n \to \infty} z_n = \overline{z} \in \R^d$ and
\[
\forall n \in \{1,2, \ldots\} \colon z_{n+1} \leq_C z_n,
\]
then 
\[
\lim_{n \to \infty} F_{X, C}\of{z_n} = F_{X, C}\of{\overline{z}}.
\]

(e) It holds
\[
\lim_{t \to \infty} F_{X, C}\of{tz} =  
\begin{cases}
    1 & \text{if} \quad z \in \qint C \\
    F_{X, C}\of{0} & \text{if} \quad z \in C\bs\qint C \\
    0 & \text{if} \quad z \not\in C  
\end{cases}
\]

\[
\lim_{t \to -\infty} F_{X, C}\of{tz} =  
\begin{cases}
    0 & \text{if} \quad z \not\in -C   \\
    F_{X, C}\of{0} & \text{if} \quad z \in -C\bs-\qint C \\
    1 & \text{if} \quad  z \in -\qint C  
\end{cases}
\]
\end{proposition}

\medskip {\sc Proof.} (a) By definition,
\begin{align*}
F_{AX + b, AC}\of{Az + b} & =  \inf_{v \in (AC)^+\bs\{0\}} P\of{v^T(AX + b) \leq v^T(Az + b)} \\
	& = \inf_{v \in (AC)^+\bs\{0\}} P\of{(A^Tv)^TX \leq (A^Tv)^Tz} \\
	& = \inf_{w \in C^+\bs\{0\}} P\of{w^TX \leq w^Tz} = F_{X, C}\of{z}
\end{align*}
since $v \in (AC)^+$ if, and only if, $v^T(Az) = (A^Tv)^Tz \geq 0$ for all $z \in C$ if, and only if, $A^Tv \in C^+$, and if $v \neq 0$ then $A^Tv \neq 0$ since $A$ is invertible and if $A^Tv \neq 0$, then $v = 0$ is not possible.

(b) If $y \leq_C z$, then $y - H^+(w) \subseteq z - H^+(w)$ for all $w \in C^+$, hence
\[
P\of{w^TX \leq w^Ty} = P\of{X \in y - H^+(w)} \leq P\of{X \in z - H^+(w)} = P\of{w^TX \leq w^Tz}.
\]
Now taking the $\inf$ over all $w \in C^+\bs\{0\}$ yields
\[
F_{X, C}\of{y} = \inf_{w \in C^+} P\of{w^TX \leq w^Ty} \leq \inf_{w \in C^+} P\of{w^TX \leq w^Tz} =F_{X, C}\of{z}.
\]

(c) If $X \leq_C Y$, then $w^TX \leq w^TY$ for all $w \in C^+$, hence
\[
F_{X,w}(z) = P(X \in z - H^+(w)) = P(w^TX \leq w^Tz) \geq P(w^TY \leq w^Tz) = F_{Y,w}(z)
\]
for all $w \in C^+$, and this implies the claimed statement.

(d) Monotonicity implies $F_{X, C}\of{z_{n+1}} \leq F_{X, C}\of{z_n}$ for all $n \in \N$, hence $\cb{F_{X, C}\of{z_n}}_{n \in \N}$ is a monotone non-increasing, bounded from below sequence, so it converges to $s := lim_{n \to \infty} F_{X, C}\of{z_n}$. Since $\overline z \leq_C z_n$ for all $n \in \N$ (because $C$ is closed), again by monotonicity $F_{X, C}\of{\overline{z}} \leq F_{X, C}\of{z_n}$ for all $n \in \N$, hence $F_{X, C}\of{\overline{z}} \leq s$. Assume that ``$<$" holds. Then, for all $w \in C^+\bs\{0\}$
\[
\forall n \in \N \colon F_{X, C}\of{\overline{z}} < s \leq F_{X, C}\of{z_n} \leq F_{X, w}\of{z_n}
\]
according to the definition of $F_{X, C}$. The definition of the infimum implies the existence of $\bar w \in C^+\bs\{0\}$ with
\[
F_{\bar w^TX}\of{\bar w^T\bar z} = F_{X, \bar w}\of{\bar z} < s \leq F_{X, \bar w}\of{z_n} = F_{\bar w^TX}\of{\bar w^Tz_n}. 
\]
The function $t \mapsto F_{\bar w^TX}\of{t}$ is right continuous and the sequence defined by $t_n := \bar w^Tz_n$ is non-increasing and convergent to $\bar w^T\bar z$. Hence
\[
F_{\bar w^TX}\of{\bar w^T\bar z} = \lim_{n \to \infty} F_{\bar w^TX}\of{\bar w^Tz_n}
\]
which produces a contradiction.

(e) Straightforward from the definitions of $F_{X, C}$, the dual cone and the quasi-interior of $C$. \pend

\medskip The importance of equivariance properties is highlighted in \cite{Serfling10JNS}. Property (a) in Proposition \ref{PropCDistribution} means that the equivariance property also involves the cone $C$: If all data points are transformed, then the ordering cone has to be transformed in the same way.

\section{Quantile functions associated to a cone}

What is a quantile for a multivariate random variable? In this section, a novel answer to this question is proposed which produces two set-valued functions as analogues to the univariate lower and upper quantile function.

If $X$ is a univariate random variable and $p \in (0,1)$, the lower $p$-quantile of $X$ is the infimum of the set
\[
\cb{r \in \R \mid P(X \leq r) \geq p},
\]
which is ``directed upward", i.e. if $r$ belongs to this set, so does $r + s$ for $s \geq 0$. On the other hand, the set 
\[
\cb{r \in \R \mid P(X < r) \leq p}
\]
is ``directed downward." The intersection of these two sets is the set of $p$-quantiles for the univariate random variable $X$. This motivates the following definitions. The symbol $\mathcal P(\R^d)$ stands for the power set of $\R^d$, i.e. the set of all subsets of $\R^d$ including $\emptyset$.

\begin{definition}
\label{DefSetLowerQuantiles} For $w \in C^+\bs\{0\}$, the function $Q^-_{X, w} \colon [0,1] \to \mathcal P(\R^d)$ defined by
\[
Q^-_{X, w}\of{p} = \cb{z \in \R^d \mid  F_{X,w}(z) \geq p}
\]
is called the {\em lower $w$-quantile function} of $X$. The function $Q^-_{X, C} \colon [0,1] \to \mathcal P(\R^d)$ defined by
\[
Q^-_{X,C}\of{p} = \cb{z \in \R^d \mid  F_{X,C}(z) \geq p}
\]
is called the {\em lower $C$-quantile function} of $X$.
\end{definition}

\begin{definition}
\label{DefSetUpperQuantiles} For $w \in C^+\bs\{0\}$, the function $Q^+_{X, w} \colon [0,1] \to \mathcal P(\R^d)$ defined by
\begin{align*}
Q^+_{X, w}\of{p} = & \cb{z \in \R^d \mid P(X \in z - \Int H^+(w) \leq p} = \\ 
= & \cb{z \in \R^d \mid P(w^TX < w^Tz) \leq p}
\end{align*}
is called the {\em upper $w$-quantile function} of $X$. The function $Q^+_{X, C} \colon [0,1] \to \mathcal P(\R^d)$ defined by
\[
Q^+_{X,C}\of{p} = \bigcap_{w \in C^+\bs\{0\}} Q^+_{X, w}\of{p}
\]
is called the {\em upper $C$-quantile function} of $X$.
\end{definition}

The following facts are immediate. First, for the sake of future reference, we formally state a simple result which makes the definitions of lower and upper quantiles completely analogous.

\begin{proposition}
\label{PropWQuantileIntersection}
It holds 
\[
\forall p \in [0,1] \colon Q^-_{X,C}\of{p} = \bigcap_{w \in C^+\bs\{0\}}Q^-_{X, w}\of{p}.
\]
\end{proposition}

{\sc Proof.} The formula follows from
\[
F_{X,C}(z) \geq p \quad \Leftrightarrow \quad \forall w \in C^+\bs\{0\} \colon F_{X, w}(z) \geq p
\]
and the definitions of $Q^-_{X, w}$, $Q^-_{X,C}$. \pend

Secondly, 
\[
Q^+_{X,C}\of{p} = \cb{z \in \R^d \mid \sup_{w \in C^+\bs\{0\}} P(w^TX < w^Tz) \leq p}.
\]
Because of this, $H^+(w) = -H^+(-w)$ and
\[
P(X \in z - \Int H^+(w)) = 1 - P(X \in z + H^+(w)),
\]
one may conclude
\begin{align*}
Q^+_{X,C}\of{p} & = \cb{z \in \R^d \mid \sup_{w \in C^+\bs\{0\}} 1 - P(X \in z + H^+(w)) \leq p} \\
		& = \cb{z \in \R^d \mid 1 - \inf_{w \in C^+\bs\{0\}} P(X \in z + H^+(w)) \leq p} \\
		& = \cb{z \in \R^d \mid F_{X, -C}(z) \geq 1 - p} = Q^-_{X,-C}\of{1-p}.
\end{align*}

\begin{remark}
\label{RemQPlusMinus}
The last formula, i.e. 
\begin{equation}
\label{EqPosNegQuantile}
\forall p \in (0,1) \colon Q^+_{X,C}\of{p} = Q^-_{X,-C}\of{1-p}
\end{equation}
means that results for lower quantiles can easily be transferred into results for upper quantiles by observing that $C$ has to be replaced by $-C$ (hence any $w \in C^+\bs\{0\}$ by $-w$) and $p$ by $1-p$. Below, we will frequently make use of this procedure. 
\end{remark}

\begin{remark}
\label{RemSurvival}
If $X$ is univariate, then
\[
\cb{r \in \R \mid F_X\of{r} \geq p} = \cb{r \in \R \mid \bar F_X\of{r} \leq 1- p}
\]
where $\bar F_X\of{r} = 1 - F_X\of{r} = P(X > r)$ is the survival function associated to $X$. In the multivariate case, a naive definition of quantiles via the joint distribution and its (joint) survival function function leads to two different concepts, see e.g. \cite{EmbrechtsPuccetti06JMA, CousinDiBernardino13JMA}, usually called the lower and the upper orthant quantile, value at risk, stochastic order etc. Our approach provides a remedy for this dilemma. Indeed, 
\begin{align*}
Q^-_{X, C}(p) & = \cb{z \in \R^d \mid \inf_{w \in C^+\bs\{0\}} F_{X, w}(z) \geq p} \\
	& = \cb{z \in \R^d \mid \inf_{w \in C^+\bs\{0\}}\sqb{1 -  \bar F_{w^TX}(w^Tz)} \geq p} \\
	& = \cb{z \in \R^d \mid 1 - \sup_{w \in C^+\bs\{0\}}\bar F_{w^TX}(w^Tz) \geq p} \\
	& = \cb{z \in \R^d \mid \sup_{w \in C^+\bs\{0\}}\bar F_{w^TX}(w^Tz) \leq 1 - p},
\end{align*}
i.e. it is also possible to generated $Q^-_{X, C}$ via survival functions, and this does not produce ambiguity. We call the function $z \mapsto \sup_{w \in C^+\bs\{0\}}\bar F_{w^TX}(w^Tz)$ the cone survival function of $X$ (with respect to the cone $C$).
\end{remark}

According to the Remark \ref{RemQPlusMinus} it is sufficient to study lower quantile functions which is done in the following. The next result shows that $Q^-_{X, w}$ is halfspace-valued.

\begin{proposition}
\label{PropWQuantile}
The function $p \mapsto Q^-_{X, w}\of{p}$ has closed convex values. In particular, 
\[
\forall p \in [0,1] \colon Q^-_{X, w}\of{p} \oplus H^+(w) = Q^-_{X, w}\of{p}.
\]
\end{proposition}

{\sc Proof.} The function $F_{X,w}$ is the composition of the non-decreasing, right continuous distribution function of $w^TX$ and the linear function $z \mapsto w^Tz$, so it is upper semicontinuous and hence has closed upper level sets which means that $Q^-_{X, w}\of{p}$ is closed.

Take $z \in Q^-_{X, w}\of{p}$ and $y \in H^+(w)$. Then $y - H^+(w) \supseteq -H^+(w)$ and hence
\[
F_{X,w}(z + y) = P(X  \in z + y - H^+(w)) \geq P(X  \in z - H^+(w)) = F_{X,w}(z)
\]
and therefore, $z+y \in Q^-_{X, w}\of{p}$. \pend

\begin{proposition}
\label{PropCQuantile} The lower $C$-quantile function $Q^-_{X, C}$ has the following properties:

(a) The function $p \mapsto Q^-_{X, C}\of{p}$ has closed convex values and satisfies 
\[
\forall p \in [0,1] \colon Q^-_{X,C}\of{p} \oplus C = Q^-_{X,C}\of{p}.
\]
In particular, $Q^-_{X, C}\of{p}$ is a connected set for each $p \in [0,1]$.

(b) For all $b \in \R^d$ and all invertible matrices $A \in \R^{d \times d}$ it holds 
\[
\forall p \in [0, 1] \colon Q^-_{AX+b, AC}\of{p} = AQ^-_{X, C}\of{p} + b.
\]

(c) If $p_1, p_2 \in [0, 1]$, $p_1 \geq p_2$, then $Q^-_{X, C}\of{p_1} \subseteq Q^-_{X, C}\of{p_2}$.

(d) If $X \leq_C Y$, then $Q^-_{X, C}\of{p} \supseteq Q^-_{Y, C}\of{p}$ for all $p \in [0,1]$.

\end{proposition}

{\sc Proof.} (a) Since $Q^-_{X, C}\of{p}$ is the intersection of closed halfspaces (see Proposition \ref{PropWQuantileIntersection} and Proposition \ref{PropWQuantile}), it is closed and convex. The formula can be proven in a similar way as the corresponding formula in Proposition \ref{PropWQuantile}. 

(b) Using the definitions and Proposition \ref{PropCDistribution} (a) we obtain
\begin{align*}
Q^-_{AX+b, AC}\of{p} 
	& = \cb{z \in \R^d \mid  F_{AX+b,AC}(z) \geq p} \\
	& = \cb{z - b \in \R^d \mid F_{AX+b,AC}(z - b + b) \geq p} + b \\
	& = \cb{z \in \R^d \mid F_{AX+b,AC}(z + b) \geq p} + b \\
	& = A\cb{A^{-1}z \in \R^d \mid F_{AX+b,AC}(A(A^{-1}z) + b) \geq p} + b \\
	& = A\cb{y \in \R^d \mid  F_{X,C}(y) \geq p} + b = AQ^-_{X, C}\of{p} + b.
\end{align*}

(c) The proof is immediate from Proposition \ref{PropCDistribution} (b). 

(d) This follows from the definition of the lower $C$-quantile and Proposition \ref{PropCDistribution} (c). \pend 
 
Proposition \ref{PropCQuantile} (a) means that the lower quantile function $p \mapsto Q^-_{X, C}\of{p}$ actually maps into the complete lattice $(\G(\R^d, C), \supseteq)$ (see \cite{HamelEtAl15Incoll} or Section \ref{SecSOP} below for definitions).  It is precisely this fact that admits to handle set-valued quantiles in the same way as scalar quantiles for univariate random variables. The parallel result for upper quantiles reads as follows.

\begin{proposition}
\label{PropUpperQuantiles}
(a) The function $p \mapsto Q^+_{X, w}\of{p}$ has convex values. In particular, 
\[
\forall p \in [0,1] \colon Q^+_{X, w}\of{p} \ominus H^+(w) = Q^+_{X, w}\of{p}.
\]

(b) The function $p \mapsto Q^+_{X, C}\of{p}$ has convex values. In particular, 
\[
\forall p \in [0,1] \colon Q^+_{X,C}\of{p} \ominus C = Q^+_{X,C}\of{p}.
\]
\end{proposition}

{\sc Proof.} (a) Follows from Remark \ref{RemQPlusMinus} and Proposition \ref{PropWQuantile}, Proposition \ref{PropCQuantile} (a).\pend

This result means that the upper quantile function $p \mapsto Q^+_{X, C}\of{p}$ actually maps into $(\G(\R^d, -C), \subseteq)$. The reader may now easily transfer the remaining properties of Proposition \ref{PropCQuantile} into ones for upper quantiles.

For a univariate $X$, the set of $p$-quantiles is 
\[
\cb{r \in \R \mid P(X \leq r) \geq p} \cap \cb{r \in \R \mid P(X < r) \geq p},
\]
and it might be tempting to define quantiles of multivariate $X$ by taking $Q^-_{X, C}\of{p} \cap Q^+_{X, C}\of{p}$. It turns out that this is asking too much: In contrast to the univariate case, this intersection can be empty which again is a consequence of the non-totalness of the order generated by $C$. 

\begin{example}
\label{ExFourPointUni}
Let $\Omega = \cb{(-1,2)^T, (0,0)^T, (1,1)^T, (2,-1)^T} \subseteq \R^2$ with a uniform distribution over $\Omega$ and $C = \R^2_+$. The following pictures show the lower and upper quantile sets for the seven cases $p \in (0, \frac{1}{4})$, $p = \frac{1}{4}$, $p \in (\frac{1}{4}, \frac{1}{2})$, $p = \frac{1}{2}$, $p \in (\frac{1}{2}, \frac{3}{4})$, $p = \frac{3}{4}$ and $p \in (\frac{3}{4}, 1)$. The intersection $Q^-_{X, C}\of{p} \cap Q^+_{X, C}\of{p}$ is non-empty exactly at the ``borderline" cases $p \in \cb{\frac{1}{4}, \frac{1}{2}, \frac{3}{4}}$. 

\begin{figure}[H]
\centering
\begin{minipage}{.45\textwidth}
  \centering
  \includegraphics[width=0.98\linewidth]{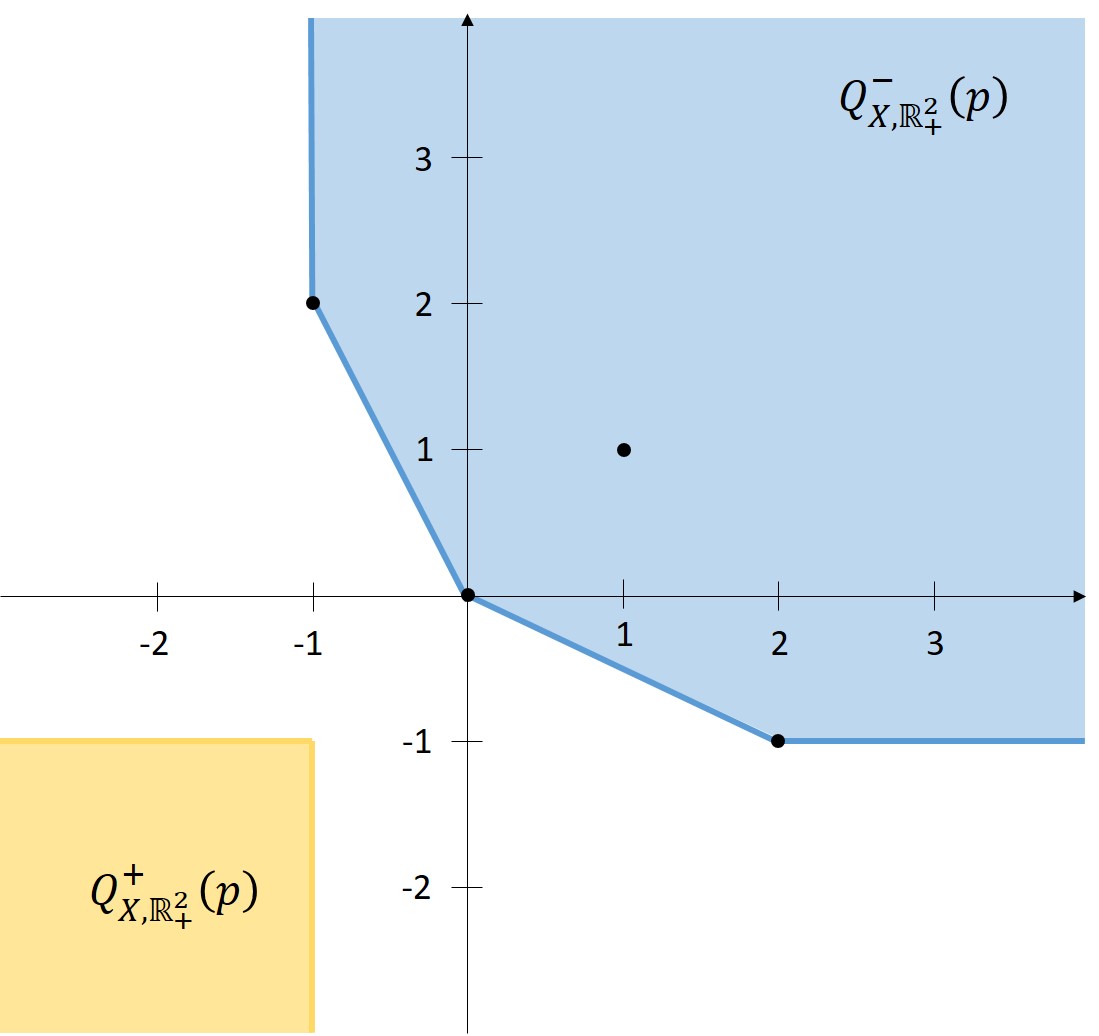}
  \caption{$p \in (0,\frac{1}{4})$}
  \label{fig:ExFour01}
\end{minipage}
\begin{minipage}{.45\textwidth}
  \centering
  \includegraphics[width=0.98\linewidth]{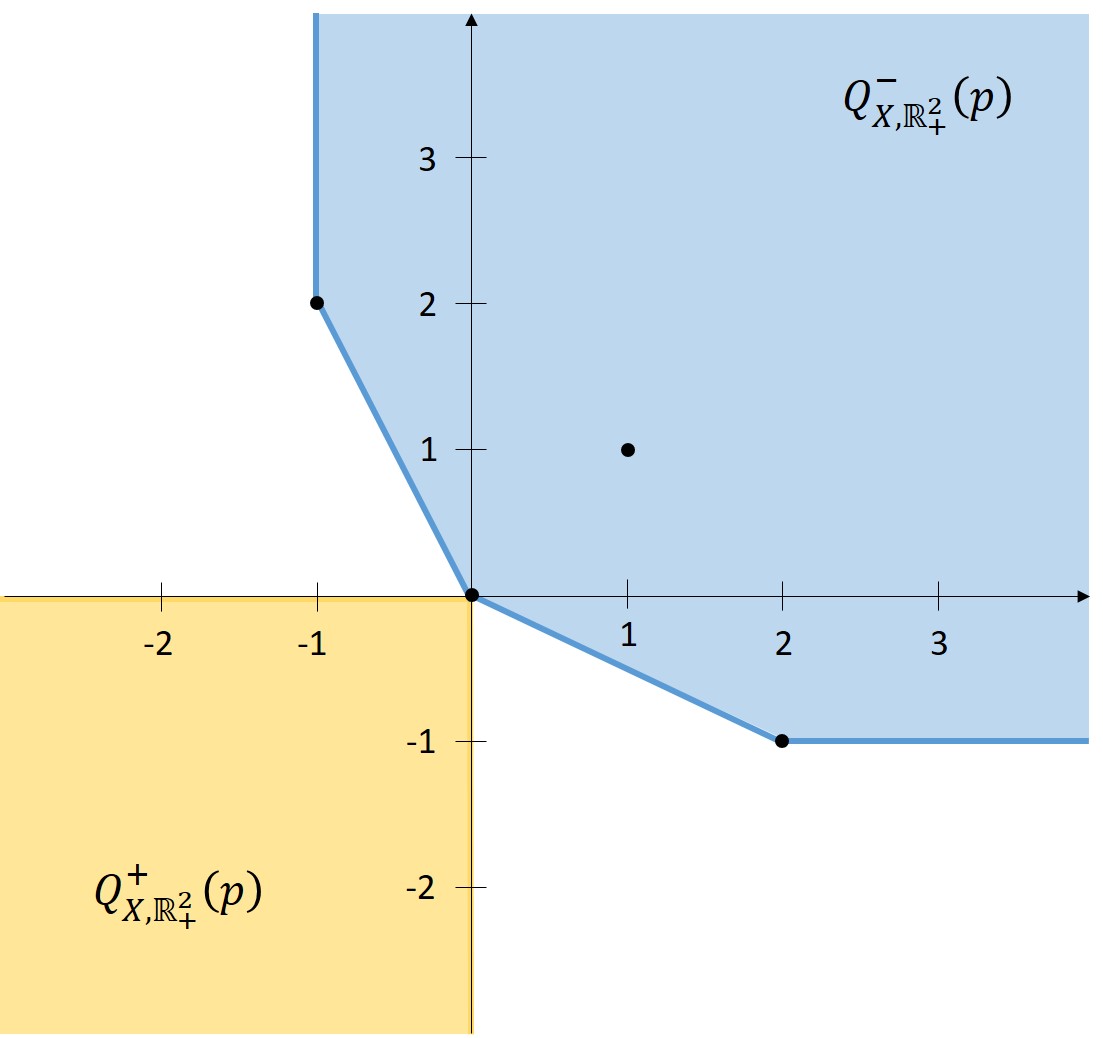}
  \caption{$p = \frac{1}{4}$}
  \label{fig:ExFour1}
\end{minipage}
\end{figure}

\begin{figure}[H]
\centering
\begin{minipage}{.45\textwidth}
  \centering
  \includegraphics[width=0.98\linewidth]{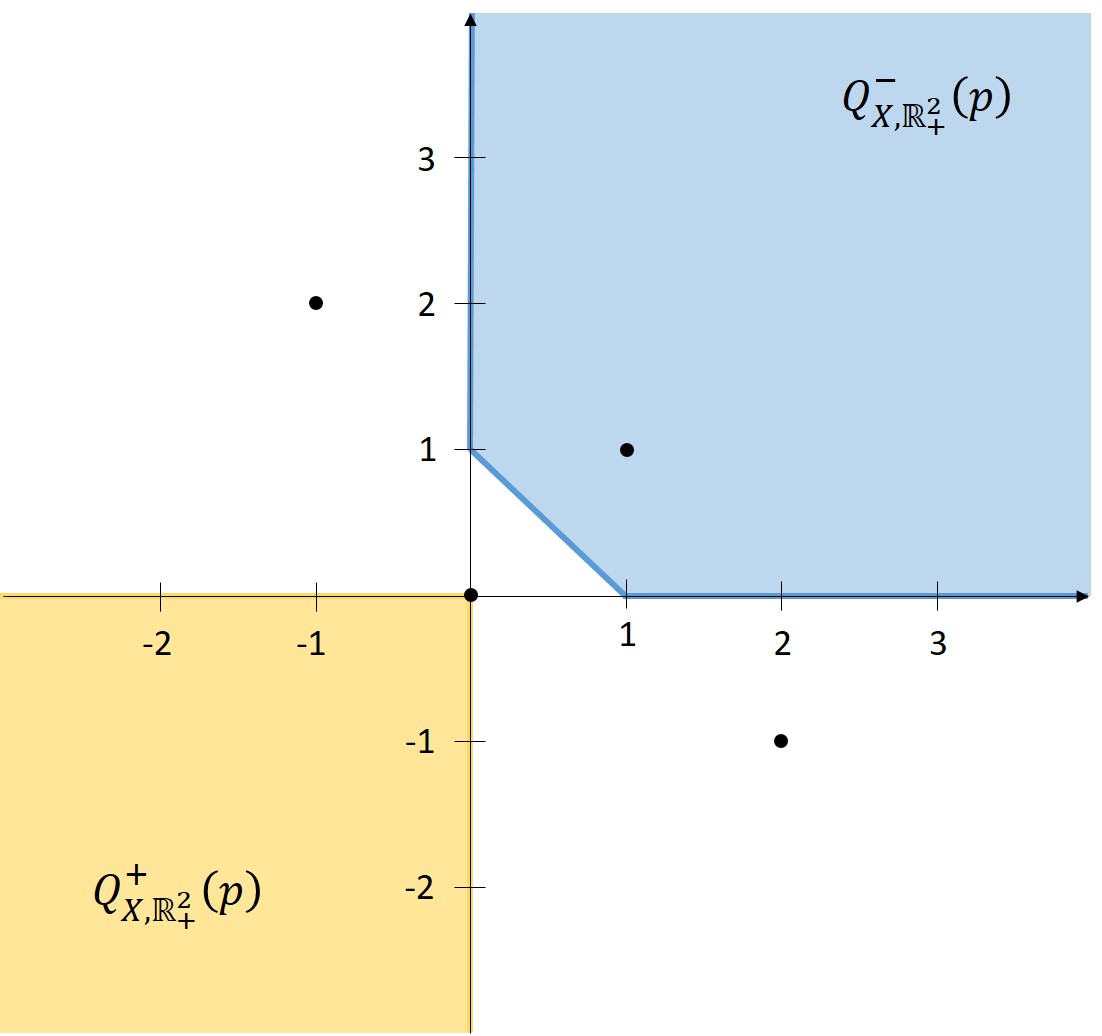}
  \caption{$p \in (\frac{1}{4},\frac{1}{2})$}
  \label{fig:ExFour12}
\end{minipage}
\begin{minipage}{.45\textwidth}
  \centering
  \includegraphics[width=0.98\linewidth]{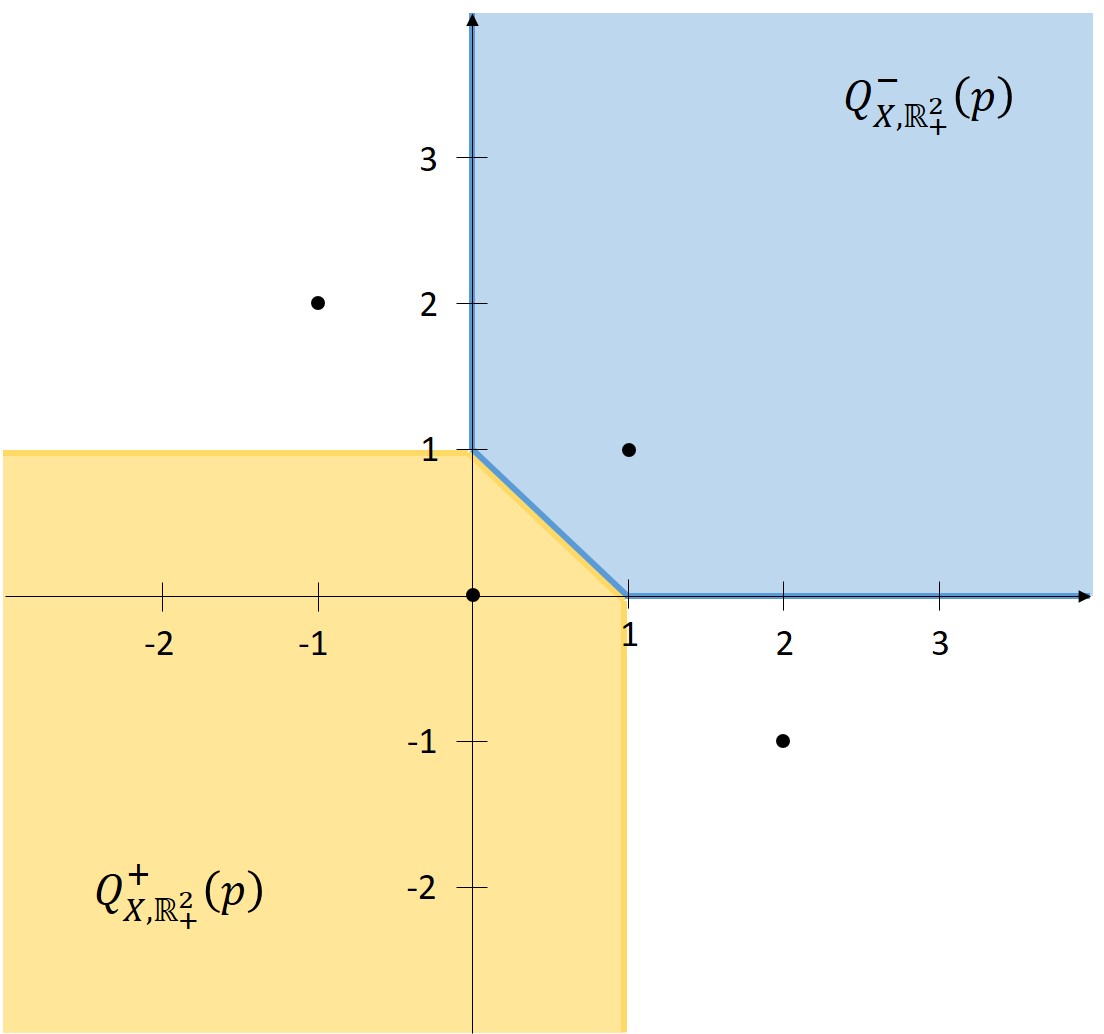}
  \caption{$p = \frac{1}{2}$}
  \label{fig:ExFour2}
\end{minipage}
\end{figure}

\begin{figure}[H]
\centering
\begin{minipage}{.45\textwidth}
  \centering
  \includegraphics[width=0.98\linewidth]{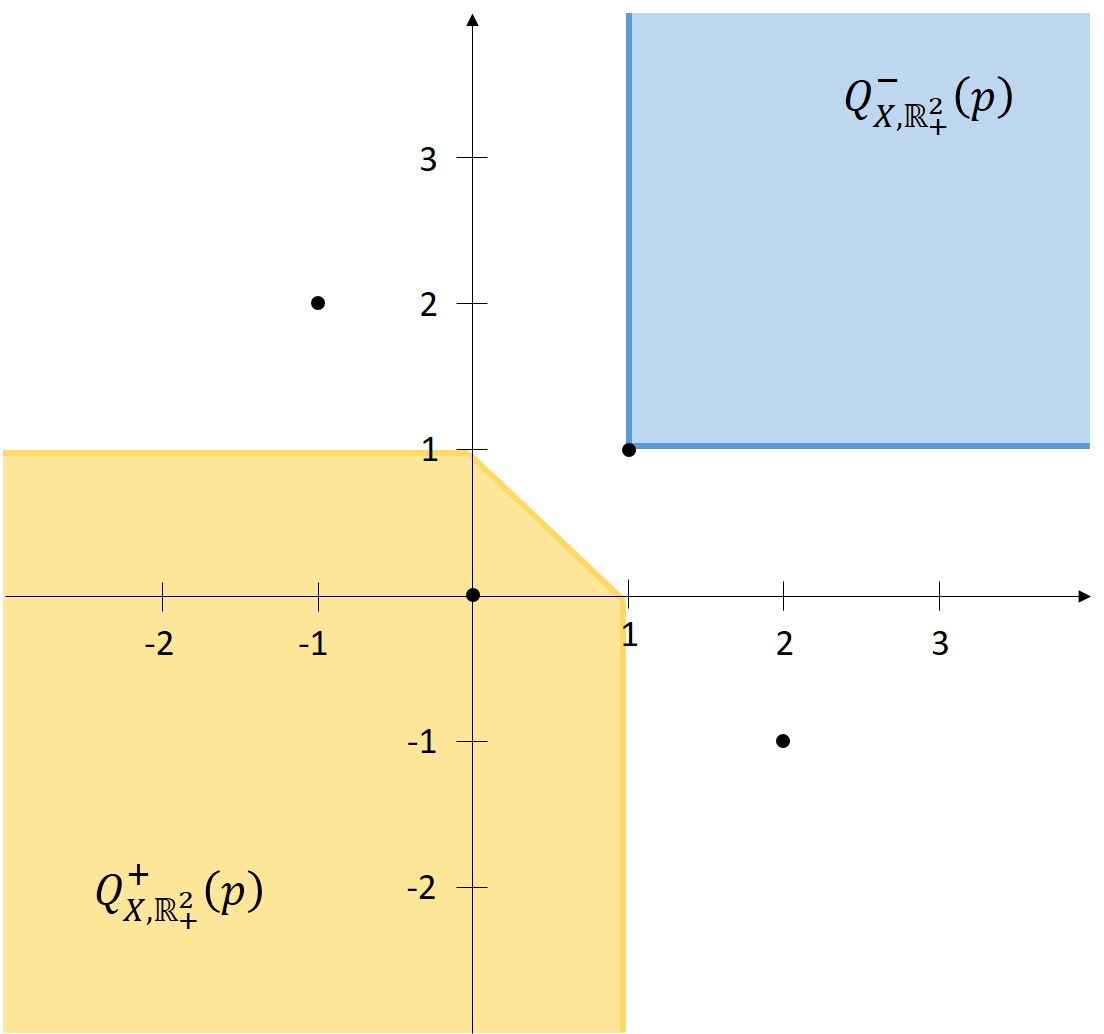}
  \caption{$p \in (\frac{1}{2},\frac{3}{4})$}
  \label{fig:ExFour23}
\end{minipage}
\begin{minipage}{.45\textwidth}
  \centering
  \includegraphics[width=0.98\linewidth]{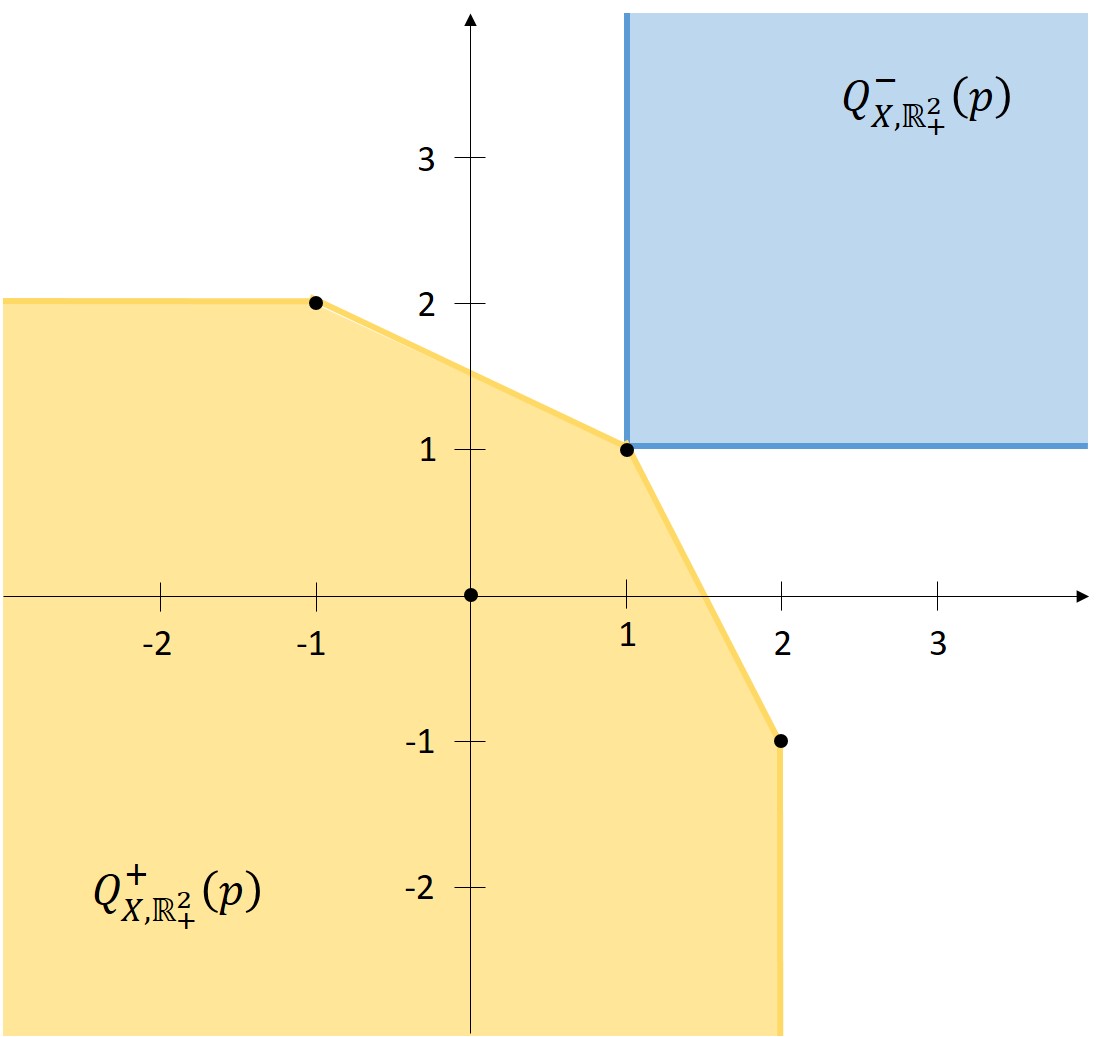}
  \caption{$p = \frac{3}{4}$}
  \label{fig:ExFour3}
\end{minipage}
\end{figure}

\begin{figure}[H]
  \centering
  \includegraphics[width=0.45\linewidth]{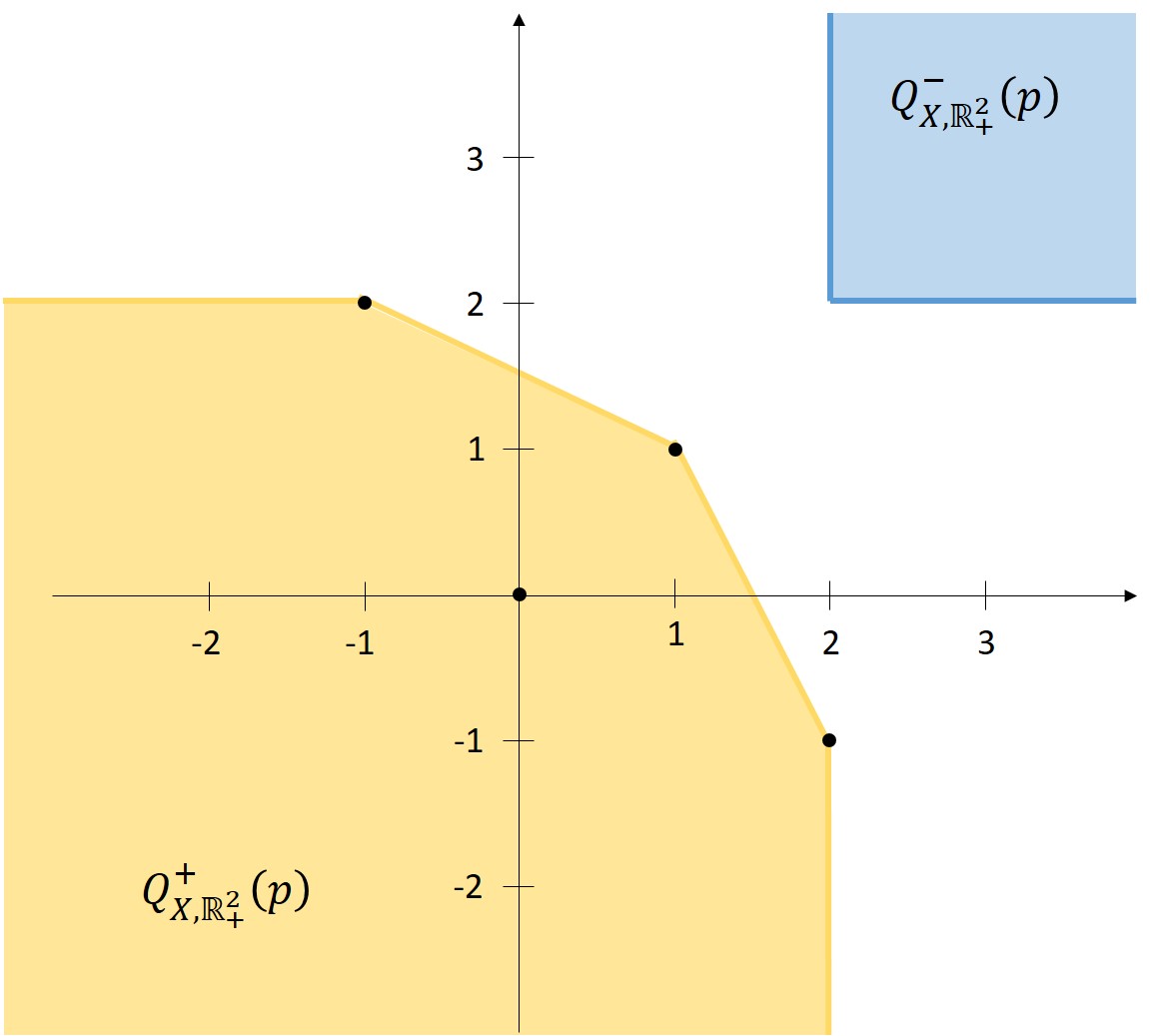}
  \caption{$p \in (\frac{3}{4},1)$}
  \label{fig:ExFour34}
\end{figure}
\end{example} 

There is dual way of writing $Q^-_{X,C}$. The proof is prepared by the following lemma which should be known (and is implicitly part of the proof of Theorem 2.11 in \cite{ZuoSerfling00AS}). The result itself is inspired by \cite[Theorem 4.1]{HallinEtAl10AS}, but we do not need the uniqueness assumption imposed therein.

\begin{lemma}
\label{LemWDistDual}
For all $w \in \R^d\bs\{0\}$ and all $z \in \R^d$ with $P(X \in z - H^+(w)) < p$ there is $y \in z + \Int H^+(w)$ such that $P(X \in y - \Int H^+(w)) < p$.
\end{lemma}

{\sc Proof.} Fix $w \in \R^d\bs\{0\}$ and $z \in \R^d$ with $P(X \in z - H^+(w)) < p$. Take $\bar z \in \R^d$ with $w^T\bar z = 1$ which exists since $w \neq 0$. Then $s\bar z \in \Int H^+(w)$ for all $s > 0$. Define $y_n = z + \frac{1}{n}\bar z \in z + \Int H^+(w)$ for $n = 1, 2, \ldots$ Then 
\[
w^Ty_n = w^T(z + \frac{1}{n}\bar z) = w^Tz + \frac{1}{n},
\]
so $w^Ty_{n+1} < w^Ty_n$ and $\lim_{n \to \infty} w^Ty_n = w^Tz$. Since $s \mapsto F_{w^TX}(s)$ is right-continuous, it follows
\[
P(X \in z - H^+(w)) = F_{w^TX}(w^Tz) = \lim_{n \to \infty}F_{w^TX}(w^Ty_n) < p,
\]
so there is $\bar n \in \cb{1, 2, \ldots}$ with
\[
F_{w^TX}(w^Ty_{\bar n}) =  P(X \in y_{\bar n} - H^+(w)) < p,
\]
hence 
\[
P(X \in y_{\bar n} - \Int H^+(w)) \leq P(X \in y_{\bar n} - H^+(w)) < p
\]
which proves the claim with $y = y_{\bar n}$. \pend

\begin{proposition}
\label{PropCQuantileDual}
For all $p \in [0,1]$,
\begin{align*}
Q^-_{X,C}\of{p} & = \bigcap_{w \in C^+\bs\{0\}}\bigcap_{y \in \R^d}\cb{y + H^+(w) \mid P\of{X \in y + H^+(w)} > 1-p} \\
	& = \bigcap_{w \in C^+\bs\{0\}}\bigcap_{y \in \R^d}\cb{y + H^+(w) \mid P\of{X \in y - \Int H^+(w)} < p}.
\end{align*}
\end{proposition}

{\sc Proof.} The two expressions on the right hand side clearly coincide since $P\of{X \in y + H^+(w)} = 1 - P\of{X \in y - \Int H^+(w)}$.

First, assume $z \not\in \bigcap_{w \in C^+\bs\{0\}}\bigcap_{y \in \R^d}\cb{y + H^+(w) \mid P\of{X \in y - \Int H^+(w)} < p}$. Then, there are $w \in C^+\bs\{0\}$, $y \in \R^d$ such that $P\of{X \in y - \Int H^+(w)} < p$ and $z \not\in y + H^+(w)$. It follows $z \in y - \Int H^+(w)$ which implies $z - H^+(w) \subseteq y - \Int H^+(w)$, so
\[
P(X \in z - H^+(w)) \leq P(X \in y - \Int H^+(w)) < p,
\]
hence $z \not\in Q^-_{X,C}\of{p}$. 

Therefore, $Q^-_{X,C}\of{p} \subseteq \bigcap\limits_{w \in C^+\bs\{0\}}\bigcap\limits_{y \in \R^d}\cb{y + H^+(w) \mid P\of{X \in y - \Int H^+(w)} < p}$.

Secondly, assume 
\[
\bar z \notin Q^-_{X,C}\of{p} = \bigcap_{w \in C^+\bs\{0\}}Q^-_{X, w}\of{p} = \bigcap_{w \in C^+\bs\{0\}} \cb{z \in \R^d \mid  F_{X,w}(z) \geq p}.
\]
Then, there is $w \in C^+\bs\{0\}$ such that $F_{X,w}(\bar z) = P(X \in \bar z - H^+(w)) < p$. Proposition \ref{LemWDistDual} yields $\bar y \in \bar z + \Int H^+(w)$ satisfying $P(X \in \bar y - \Int H^+(w)) < p$. If 
\[
\bar z \in \bigcap\limits_{w \in C^+\bs\{0\}}\bigcap\limits_{y \in \R^d}\cb{y + H^+(w) \mid P\of{X \in y - \Int H^+(w)} < p}
\]
would be true, then also $\bar z \in \bar y + H^+(w)$ and
\[
\bar z \in \of{\bar y -\Int H^+(w)} \cap \of{\bar y + H^+(w)},
\]
which is a contradiction. So, $\bar z \not\in \bigcap\limits_{w \in C^+\bs\{0\}}\bigcap\limits_{y \in \R^d}\cb{y + H^+(w) \mid P\of{X \in y - \Int H^+(w)} < p}$. This shows $Q^-_{X,C}\of{p} \supseteq \bigcap\limits_{w \in C^+\bs\{0\}}\bigcap\limits_{y \in \R^d}\cb{y + H^+(w) \mid P\of{X \in y - \Int H^+(w)} < p}$. 
\pend

The previous result can also easily be transferred into a dual representation of the upper quantiles using the device of Remark \ref{RemQPlusMinus}.

In the following example, the quantile sets $Q_X^{jdf}\of{p}$ based on the joint distribution are compared to the lower $C$-quantiles. There are two important insights. First, $Q_X^{jdf}\of{p}$ is in general non-convex. Second, the lower $C$-quantiles are not even equal to the convex hulls of the joint distribution quantiles in general.

\begin{example}
\label{Joint}
Consider the four-point uniform distribution of Example \ref{ExFourPointUni}. The figures \ref{fig:ExFourJointC1}, \ref{fig:ExFourJointC2} compare the lower $C$-quantile with the joint distribution quantile defined by
\[
Q_X^{jdf}\of{p}=\cb{z \in \R^2 \mid F^{jdf}_{X}(z) \geq p},
\]
namely  $Q^-_{X,C}\of{\frac{1}{4}}$ versus $Q_X^{jdf}\of{\frac{1}{4}}$ and $Q^-_{X,C}\of{\frac{1}{2}}$ versus $Q_X^{jdf}\of{\frac{1}{2}}$, respectively. It can be seen that the joint distribution quantiles are non-convex and that the corresponding lower $C$-quantile sets are not their convex hulls in general. Note that the ``blue" sets also cover the ``green" ones.

\begin{figure}[H]
\centering
\begin{minipage}{.45\textwidth}
  \centering
  \includegraphics[width=0.98\linewidth]{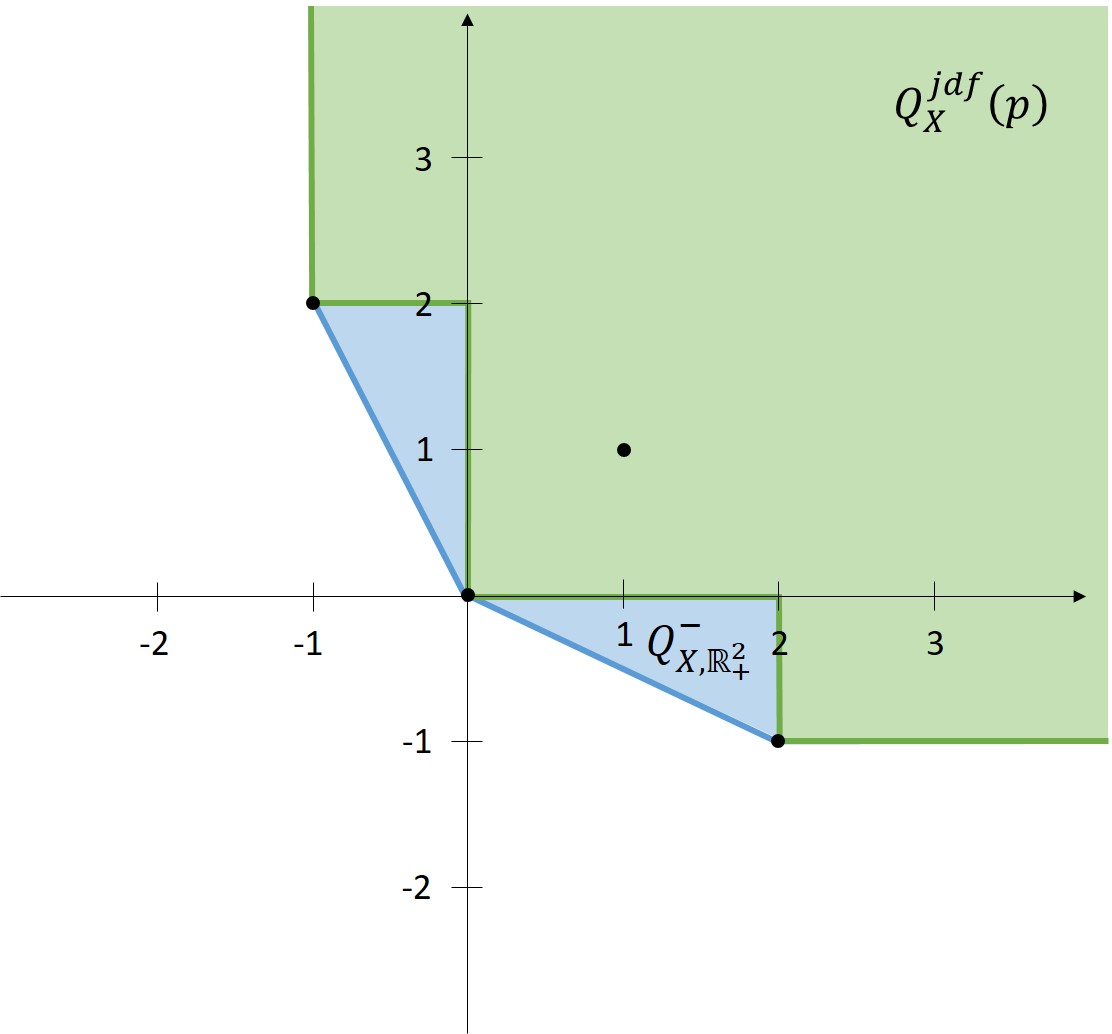}
  \caption{$p = \frac{1}{4}$}
  \label{fig:ExFourJointC1}
\end{minipage}
\begin{minipage}{.45\textwidth}
  \centering
  \includegraphics[width=0.98\linewidth]{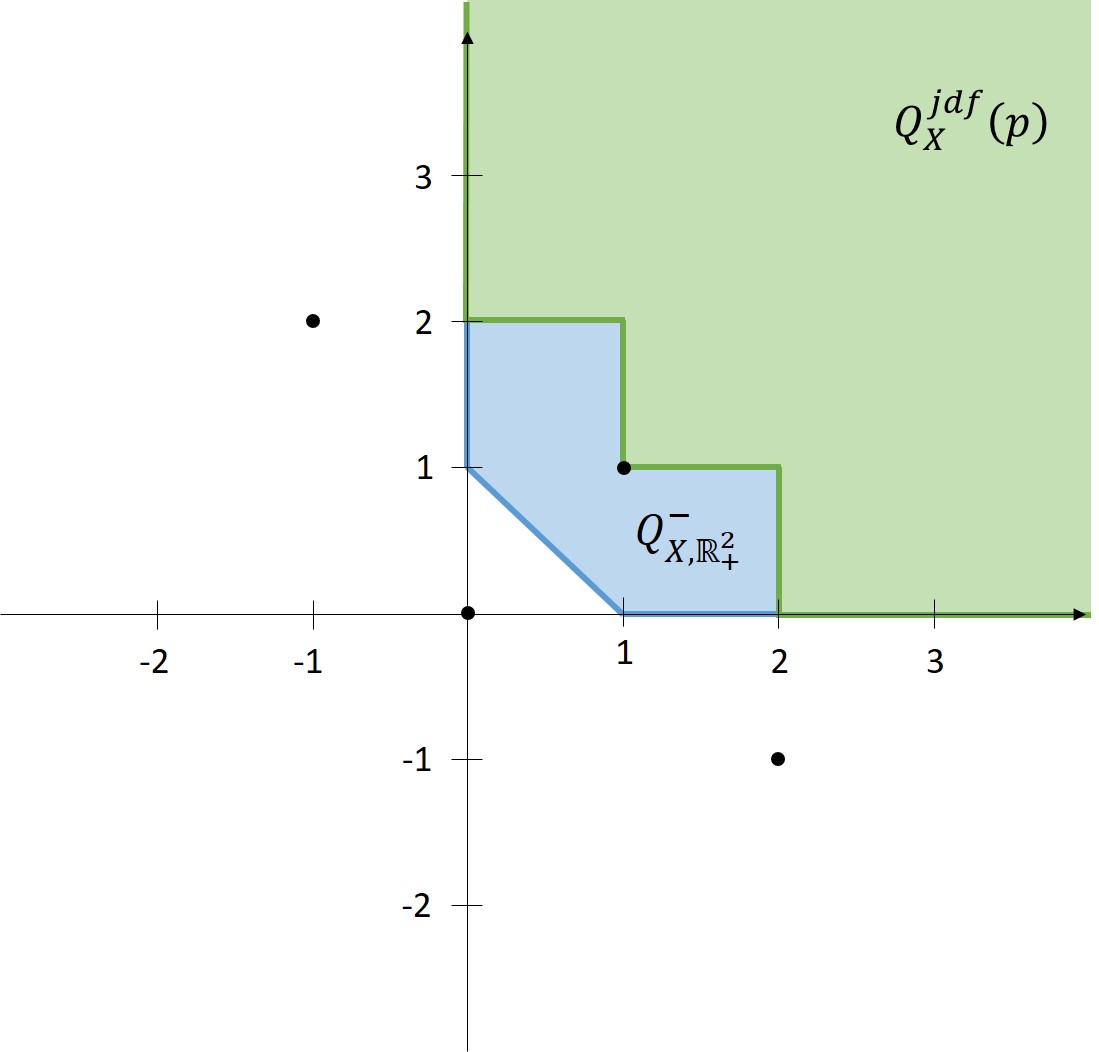}
  \caption{$p = \frac{1}{2}$}
  \label{fig:ExFourJointC2}
\end{minipage}
\end{figure}
\end{example}

The following example illustrates again that our quantiles are different from the ones defined via the joint distribution function. Moreover, the lower $C$-quantiles are also different from a component-wise defined quantile even if $C = \R^d_+$, i.e. different from the set $(q_{X_1}\of{p}, \ldots,  q_{X_d}\of{p})^T + \R^d_+$ where $q_{X_i}\of{p}$ is the univariate lower quantile of $X_i$ for $i \in \cb{1, \ldots, d}$.

\begin{example}
\label{ExChiSquare}
The following pictures show different quantiles for $\alpha = 0.5$ for the bivariate non-central $\chi^2$-distribution with $k=1$ and a non-centrality parameter close to 0 with $C = \R^2_+$. The {\color{red} red} set represents the component-wise quantile, the {\color{blue} blue} set the lower $\R^d_+$-quantile and the {\color{green} green} set is $Q_X^{jdf}\of{0.5}$. The bigger sets cover the smaller ones. Figure \ref{fig:chi} shows that the component-wise quantile is a superset of the lower $C$-quantile, whereas the lower $C$-quantile is a superset of the joint distribution quantile. 

\begin{figure}[H]
  	\centering
  	\includegraphics[width=0.7\linewidth]{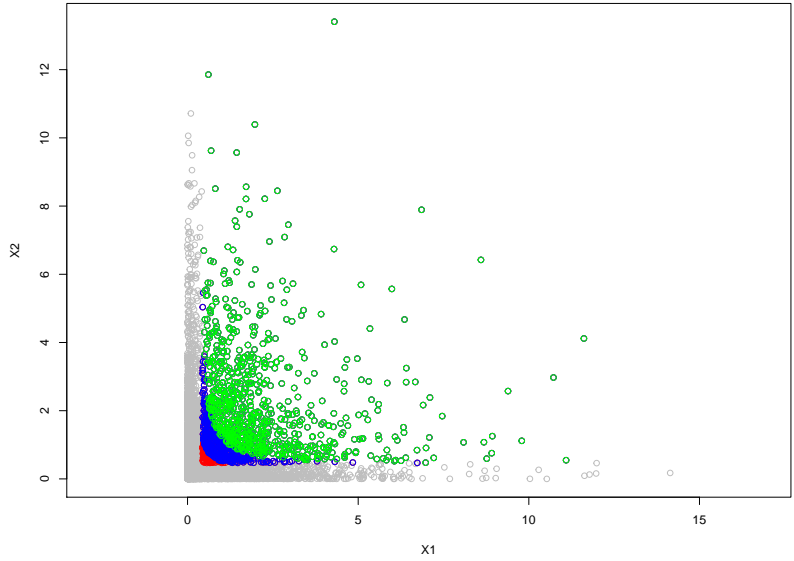}
  	\caption{Bivariate $\chi^{2}$-Distribution}
  	\label{fig:chi}
\end{figure}

However, the component-wise quantile and lower $C$-quantile may coincide for some distributions, as shown by figure \ref{fig:Norm}.

\begin{figure}[H]
	\centering
  	\includegraphics[width=.7\linewidth]{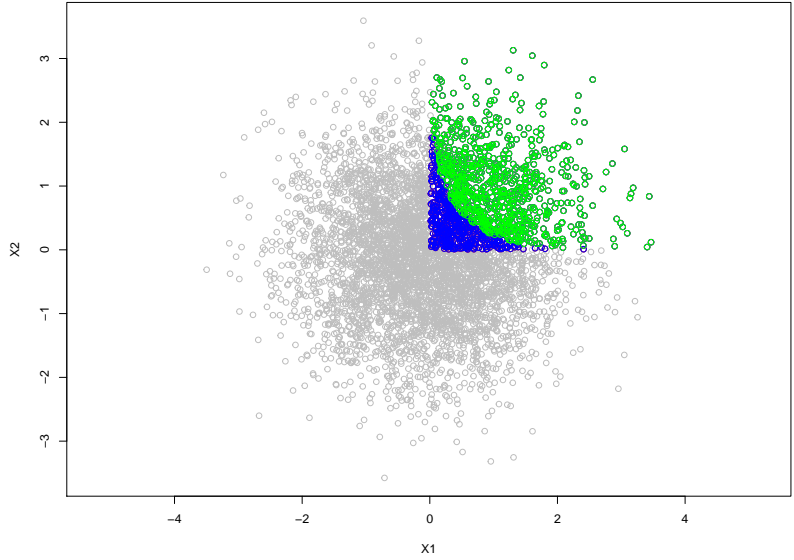}
  	\caption{Bivariate Standard Normal Distribution, $p = 0.5$}
  	\label{fig:Norm}
\end{figure}
\end{example}

\section{Relationships with univariate quantiles}

\label{SecScalarQuantile}

The lower quantile of the (univariate) random variable $w^TX$ for $w \in C^+\bs\{0\}$ and $X \colon \Omega \to \R^d$ is the function $q^-_{w^TX} \colon [0,1] \to \R$ defined by
\[
q^-_{w^TX}\of{p} = \inf \cb{s \in \R \mid P\of{w^TX \leq s} \geq p}.
\]
The lower $C$-quantile of $X$ can also be expressed in terms of the family $\cb{q^-_{w^TX}}_{w \in C^+\bs\{0\}}$. The result reads as follows.
 
\begin{proposition}
\label{PropQuantileScalar}
(a) If $w \in \R^d\bs\{0\}$, then
\[
\forall p \in [0,1] \colon Q_{X,w}(p) = \cb{z \in \R^d \mid w^Tz \geq q^-_{w^TX}(p)}.
\]

(b) It holds
\[
\forall p \in [0,1] \colon Q^-_{X,C}(p) = \bigcap_{w \in C^+\bs\{0\}} \cb{z \in \R^d \mid w^Tz \geq q^-_{w^TX}(p)}.
\]

(c) If $w \in C^+\bs\{0\}$, then
\[
\forall p \in [0,1] \colon q^-_{w^TX}(p) \leq \inf_{z \in Q^-_{X,C}(p)}w^Tz.
\]
\end{proposition}

{\sc Proof.} (a) If $z \in Q_{X,w}\of{p}$, then $P(w^TX \leq w^Tz) \geq p$, hence $q^-_{w^TX}(p) \leq w^Tz$ by definition of $q^-_{w^TX}$.

Conversely, if $q^-_{w^TX}(p) \leq w^Tz$ for $z \in \R^d$, then by monotonicity of $F_{w^TX}$ and a known property of quantile functions
\[
F_{w^TX}(w^Tz) \geq F_{w^TX}(q^-_{w^TX}(p)) \geq p,
\]
hence $z \in Q_{X,w}(p)$.

(b) This follows immediately from (a) and Proposition \ref{PropWQuantileIntersection}. 

(c) Is immediate from (b). \pend 

Loosely speaking, (b) of the preceding proposition means that $Q^-_{X,C}$ can be constructed from the family $\cb{q^-_{w^TX}}_{w \in C^+}$ of scalar quantiles which provides a toehold for a computational approach. On the other hand, the function 
\[
w \mapsto \inf_{z \in Q^-_{X,C}(p)}w^Tz
\]
coincides up to signs with the (sublinear) support function of the closed convex set $Q^-_{X,C}(p)$ which yields
\[
\forall p \in [0,1] \colon  Q^-_{X,C}\of{p} = \bigcap_{w \in C^+\bs\{0\}}\cb{z \in \R^d \mid w^Tz \geq \inf_{z \in Q^-_{X,C}(p)}w^Tz}. 
\]
The latter formula and the one in (b) produce the same set by different ``scalarization functions." Observe that the functions
$w \mapsto \inf_{z \in Q^-_{X,C}(p)}w^Tz$ are superlinear while the functions $w \mapsto q^-_{w^TX}(p)$ are not. An example for strict inequality in (c) is as follows.

\begin{example}
Take the 4-point uniform distribution from Example \ref{ExFourPointUni}, $p = 3/8$, $w = (2,1)^T$. Then, on the one hand, $q^-_{w^TX}(\frac{3}{8}) = 0$ and $\inf_{z \in Q^-_{X,C}(\frac{3}{8})}w^Tz = 1$. On the other hand, for $w = (1,1)^T$ we obtain $q^-_{w^TX}(\frac{3}{8}) = 1 = \inf_{z \in Q^-_{X,C}(\frac{3}{8})}w^Tz$.
\end{example}

As far as the relationship to the scalar upper quantiles is concerned we have the following result.

\begin{proposition}
\label{PropUpperWquantileScalar}
For all $p \in (0,1)$,
\[
\forall w \in C^+\bs\{0\} \colon 
	Q^+_{X,w}(p ) =  \cb{z \in \R^d \mid w^Tz \leq q^+_{w^TX}(p )} 
\]
where $q^+_{w^TX}(p )= \sup\cb{s \in \R \mid P\of{w^TX < s} \leq p}$. Moreover,
\[
Q^+_{X,C}(p ) = \bigcap_{w \in C^+\bs\{0\}}\cb{z \in \R^d \mid w^Tz \leq q^+_{w^TX}(p )}.
\]
\end{proposition}

{\sc Proof.} This follows from Proposition \ref{PropQuantileScalar}, Remark \ref{RemQPlusMinus} and $q^-_{-w^TX}(1-p) = -q^+_{w^TX}(p)$. \pend

\section{The link to Tukey depth}
\label{SecConnectTukey}

The Tukey halfspace depth function $HD_X \colon \R^d \to [0,1]$ associated with the random variable $X \colon \Omega \to \R^d$ defined in \cite{Tukey75} is given by
\[
HD_X(z)  = \inf_{w \in \R^d\bs\{0\}}P\of{X \in z - H^+(w)}.
\]
Clearly, the set $\R^d\bs\{0\}$ can be replaced by a unit sphere or any other set $S \subseteq \R^d\bs\{0\}$ with $\cup_{t > 0} tS = \R^d\bs\{0\}$ since $H^+(w) = H^+(tw)$ whenever $t>0$, $w \in \R^d\bs\{0\}$. The Tukey depth regions are the upper level sets of the Tukey depth function:
\[
D_{X}(p ) = \cb{z \in \R^d \mid HD_{X}\of{z} \geq p}, \; p \in [0,1].
\]
One may easily recognize the Tukey depth function as a special case of the cone distribution function for the cone $C = \{0\}$ with $C^+ = \R^d$, i.e.
\[
\forall z \in \R^d \colon  F_{X,\{0\}}\of{z}  = HD_{X}(z)
\]
as well as the Tukey depth regions as set-valued quantiles, i.e.
\[
Q^-_{X, \{0\}}(p ) = \cb{z \in \R^d \mid  F_{X,\{ 0 \}}(z) \geq p}  = D_{X}(p ).
\]
This shows that the cone distribution function can be seen as a generalization of the Tukey depth function to the case of more general order relations. One may also realize that the ``dual representation" of $Q^-_{X,C}$ generalizes the dual representations of Tukey depth regions as, for example, given in \cite{HallinEtAl10AS}. 

The Tukey depth region $D_{X}(p )$ for $0 < p <1$ is a compact set if $X$ has a continuous distribution, see \cite{ZuoSerfling00AS}. Of course, such a result cannot be expected for general cones $C$. 

One more remark on the relation between Tukey depth regions and set-valued quantiles in the univariate case might clarify matters further. On the one hand, it is clear that in the univariate case, i.e. $d=1$, $D_X(p ) = Q^-_{X, \{0\}}(p )$ is not the set of (lower) $p$-quantiles in general. On the other hand, the sets $Q^-_{X, \R_+}(p )$ and $Q^+_{X, \R_+}(p )$ coincide with the set of lower and upper $p$-quantiles, respectively, and of course $F_{X, \R_+} = F_X$ is the usual cumulative distribution function. This means that the lower $C$-distribution function is a common generalization of Tukey's depth function--which can be seen as a ``measure of centrality"--and the univariate cdf. In the same way, the lower $C$-quantiles are a common generalization of Tukey's depth regions and univariate lower quantiles. It should already become apparent from this discussion that flexibility concerning the cone $C$ ($\{0\}$ for Tukey depth, $\R_+$ for univariate quantiles) is an important feature of the theory presented in this paper.

\section{Value at risk}

Using the definition of the (set-valued) quantiles given in the previous section we introduce the Value at Risk of a multivariate position completely parallel to the scalar case (compare e.g. \cite[p. 207]{FoellmerSchied11Book3rd}). It turns out that our new Value at Risk not only enjoys the same properties as the scalar one for univariate random variable, but it has an equally intuitive financial interpretation as a risk measure.

\begin{definition}
\label{DefVaR}
Let $0 < \alpha \leq 1$. The Value at Risk of $X \colon \Omega \to \R^d$ at level $\alpha$ is
\[
VaR_\alpha(X) = Q^-_{-X,C}\of{1-\alpha}.
\]
\end{definition}

\begin{proposition}
\label{PropVaR}
(a) It holds
\begin{align*}
VaR_\alpha(X) & = \cb{z \in \R^d \mid \sup_{w \in C^+\bs\{0\}}P\of{X + z \in -\Int H^+(w)} \leq \alpha} \\
	& = \bigcap_{w \in C^+\bs\{0\}} \negthickspace\negthickspace\cb{z \in \R^d \mid P\of{X + z \in -\Int H^+(w)} \leq \alpha}
\end{align*}

(b) The function $X \mapsto VaR_\alpha(X)$ maps into $\G(\R^d, C)$, is positively homogeneous and $\R^d$-translative, i.e.
\[
\forall y \in \R^d \colon VaR_\alpha(X + y \One) = VaR_\alpha(X) - y. 
\]

(c) $X \mapsto VaR_\alpha(X)$ is monotone nonincreasing with respect to $ \leq_C$, i.e. $X \leq_C Y$ implies $VaR_\alpha(X) \subseteq VaR_\alpha(Y)$.
\end{proposition}

{\sc Proof.}\\ 
(a) The obvious fact $P\of{X + z \in H^+(w)} = 1 - P\of{X + z \in -\Int H^+(w)}$ yields
\begin{align*}
VaR_\alpha(X) & = \cb{z \in \R^d \mid F^-_{-X,C}(z) \geq 1-\alpha} \\
	& = \cb{z \in \R^d \mid \inf_{w \in C^+\bs\{0\}}P\of{-X \in z - H^+(w)} \geq 1 - \alpha} \\
	& = \cb{z \in \R^d \mid \inf_{w \in C^+\bs\{0\}}P\of{X + z \in H^+(w)} \geq 1 - \alpha} \\
	& = \cb{z \in \R^d \mid \inf_{w \in C^+\bs\{0\}}\sqb{1 - P\of{X + z \in -\Int H^+(w)}} \geq 1 - \alpha} \\
	& = \cb{z \in \R^d \mid \inf_{w \in C^+\bs\{0\}}-P\of{X + z \in -\Int H^+(w)} \geq -\alpha} \\
	& = \cb{z \in \R^d \mid -\negthickspace\negthickspace\negthickspace\inf_{w \in C^+\bs\{0\}}-P\of{X + z \in -\Int H^+(w)} \leq \alpha} \\
	& = \cb{z \in \R^d \mid \sup_{w \in C^+\bs\{0\}}P\of{X + z \in -\Int H^+(w)} \leq \alpha}.
\end{align*}

(b) Everything follows from Proposition \ref{PropCQuantile} (a), (b). 

(c) This from Proposition \ref{PropCQuantile} (d). \pend

The formula in (a) of the previous proposition has a nice financial interpretation: Let us assume that $X$ denotes a future random financial position in ``physical units," i.e. $X_i$ is the number of units of asset \#$i$ in the future portfolio for $i = 1, \ldots, d$ (see \cite{HamelHeyde10SIFIN} for more explanations and references). Then, $VaR_\alpha(X)$ contains all $z \in \R^d$, i.e. all deterministic portfolios which could be deposited in a risk free manner at initial time, such that for each "vector of relative weights" $w$ the probability of bankruptcy for the merged position $X + z$ at terminal time, i.e. for the event $w^TX + w^Tz < 0$, is at most $\alpha$. 

If the $w$'s are understood as (relative) prices, then it is of course not very realistic that they do not change over time. Therefore, it is very desirable to extend the concepts introduced above to random cones $C$ and even more general random sets (see Section \ref{SecConPer} below).

Using the scalar representation formulas for the lower $C$-quantile in Proposition \ref{PropQuantileScalar} we can give similar formulas for the Value at Risk.

\begin{corollary}
\label{CorVarScalar}
(a) It holds
\[
VaR_\alpha(X) = \bigcap_{w \in C^+\bs\{0\}}\cb{z \in \R^d \mid w^Tz \geq VaR^{sca}_\alpha(w^TX)}.
\] 

(b) Conversely,
\[
VaR^{sca}_\alpha(w^TX) \leq \inf\cb{w^Tz \mid z \in VaR_\alpha(X)}.
\]
\end{corollary}

{\sc Proof.} This follows from the definition of VaR and Proposition \ref{PropQuantileScalar}. \pend

The formulas in Corollary \ref{CorVarScalar} admit to compare the Value at Risk with previously defined concepts. First, observe that if $C = \R^d_+ = C^+$, then the unit vectors $e^i$, $i = 1, \ldots, d$, are included in $C^+\bs\{0\} = \R^d_+\bs\{0\}$, hence
\begin{align*}
VaR_\alpha(X) & = \bigcap_{w \in C^+\bs\{0\}}\cb{z \in \R^d \mid w^Tz \geq VaR^{sca}_\alpha(w^TX)} \\
	& \subseteq \bigcap_{i \in \cb{1, \ldots, d}}\cb{z \in \R^d \mid z_i \geq VaR^{sca}_\alpha(X_i)}.
\end{align*}
This means that the set-valued VaR is ``more conservative" as the component-wise VaR since there are possibly less risk compensating portfolios in $VaR_\alpha(X)$ than in the component-wise VaR. This, of course, makes sense due to effects of dependencies among the components of $X$. Example \ref{ExChiSquare} above shows that even for $C = \R^d_+$ our quantile-based VaR can be different from the component-wise one. 

Finally, one might suspect that a definition via upper quantiles (see \cite[Definition 4.45]{FoellmerSchied11Book3rd}) as
\[
VaR_\alpha(X) = -Q^+_{X,C}\of{\alpha} = \bigcap_{w \in C^+\bs\{0\}}\cb{-z \in \R^d \mid P(w^TX < w^Tz) \leq \alpha}
\]
produces another version of the Value at Risk. However, it is just a little exercise (compare Remark \ref{RemSurvival}) to show that this leads to the very same set as Definition \ref{DefVaR}. Thus, as in the scalar case, lower and upper quantiles produce the same Value at Risk. 

On the other hand, one may define the VaR via the joint distributions functions. This has been done in \cite[Definition 17]{EmbrechtsPuccetti06JMA} where the `multivariate lower-orthant (LO-)Value-at-Risk' was defined as
\[
\underline{VaR}_{\alpha}\of{X} =  \bd\cb{z \in \R^d \mid F^{jdf}_{-X}(z) \geq 1- \alpha}
\]
(notation adopted to our setting). The symbol $\bd$ stands for the topological boundary. Since $F_{-X, \R^d_+}(z) \geq F^{jdf}_{-X}(z)$ for all $z \in \R^d$ (see Remark \ref{RemCDFvsJDF}) we clearly have
\[
VaR_\alpha(X) = \cb{z \in \R^d \mid F_{-X, \R^d_+}(z) \geq 1- \alpha} \supseteq  \cb{z \in \R^d \mid F^{jdf}_{-X}(z) \geq 1- \alpha}.
\]
In this sense, our Value at Risk is ``less conservative" than the LO-Value-at-Risk. In addition, the following example discloses another important feature. While our Value at Risk has convex values, the set $\cb{z \in \R^d \mid F^{jdf}_{-X}(z) \geq 1- \alpha}$ is not convex in general. This is very hard to justify: Why is a mixture of two risk compensating portfolios not risk compensating anymore? In particular, if it is ``very close" to one of the two original portfolios? It also makes it extremely difficult to build a calculus for functions like $X \mapsto \underline{VaR}_{\alpha}\of{X}$. Similar remarks can be made about the `upper-orthant (UO-)Value-at-Risk' of \cite[Definition 17]{EmbrechtsPuccetti06JMA} which is--in contrast to our Value at Risk--different from the LO-Value-at-Risk (again, compare Remark \ref{RemSurvival}). In \cite[Section 3]{HamelHeyde10SIFIN}, it is explained that the appearance of the LO- and UO-version of Value-at-Risk is a consequence of the fact that "being strictly less" is not the same as "not being greater than or equal to" with respect to a general vector order.

In more recent works such as \cite{CousinDiBernardino13JMA}, the Embrechts/Puccetti Value-at-Risk serves as a stepping stone for the construction of a vector-valued Value at Risk where according to sophisticated criteria a single point from a set like $\underline{VaR}_{\alpha}\of{X}$ is selected. Clearly, such approaches loose information on the multivariate distribution $X$; some more remarks on this can be found in Section \ref {SecSOP}. Note also that we do not require any type of `regularity' as in \cite[p. 36]{CousinDiBernardino13JMA} which means that our definitions also work well e.g. for empirical distributions. In \cite{TorreLilloLanaido15IME}, a different idea is pursued: instead of halfspaces as in Tukey's depth function, the cone $\R^d_+$ itself (more general cones are not considered) is turned and from the resulting, in general non-convex sets points are chosen. In \cite[Section 7]{CascosMolchanov07FS}, another set-valued Value at Risk appears which is even ``less conservative" than the component-wise Value at Risk.

\section{Multivariate stochastic dominance}

In analogy to the scalar case, a definition of First Order Stochastic Dominance (FSD) based on the lower $C$-distribution function is given. Moreover, it is also shown that the FSD can be expressed in terms of the lower $C$-quantile. Therefore, this type of stochastic dominance depends on the order generated by the cone $C$, it changes if $C$ changes. Previous definitions of stochastic orders involve the joint distribution and joint survival function, respectively. This approach leads to two different versions of FSD, usually called `upper orthant order' and `lower orthant order' as in Definition 3.3.1 of the standard reference \cite{MuellerStoyan02Book}.

\begin{definition}
\label{DefFSD}
The random variable $Y \colon \Omega \to \R^d$ is said to stochastically dominate the random variable $X \colon \Omega \to \R^d$, written as $Y \succeq^C_{FSD} X$, iff
\[
\forall z \in \R^d \colon F_{Y,C}\of{z} \leq F_{X,C}\of{z}.
\]
\end{definition} 

\begin{proposition}
\label{PropPropertyFSD}
For the random variables $X, Y \colon \Omega \to \R^d$, the following statements are equivalent:

(a) $Y \succeq^C_{FSD} X$,

(b) It holds
\[
\forall p \in [0,1] \colon Q^-_{Y,C}\of{p} \subseteq Q^-_{X,C}\of{p}.
\]

(c) It holds
\[
\forall \alpha \in [0,1] \colon VaR_\alpha(X) \subseteq VaR_\alpha(Y.)
\]
\end{proposition}

{\sc Proof.} From the definition of $Q^-_{X,C}\of{p} = \cb{z \in \R^d \mid  F_{X,C}(z) \geq p}$ it is immediate that (a) implies (b). Conversely, if $F_{Y,C}\of{\bar z} > F_{X,C}\of{\bar z}$ for some $\bar z \in \R^d$, then $\bar z \in Q^-_{Y,C}\of{\bar p}$, but $\bar z \not\in Q^-_{X,C}\of{\bar p}$ for $\bar p = F_{Y,C}\of{\bar z}$ contradicting (b), so (b) implies (a). The equivalence of (b) and (c) is clear from the definition of $VaR_\alpha$ as a lower $C$-quantile. \pend

Again, as for the Value at Risk, the ambiguity between `lower orthant' and `upper orthant' orders disappears based on the observation in Remark \ref{RemSurvival}: Our stochastic dominance is an intersection of univariate stochastic dominance orders generated by $F_{X,w}$, or, equivalently, by $\bar F_{X,w}$. Finally, FSD is monotone with respect to the point-wise order (see Proposition \ref{PropCDistribution} (c)):
\[
X \leq_C Y \quad \Rightarrow \quad F_{Y,C} \leq F_{X,C} \quad \Leftrightarrow \quad Y \succeq^C_{FSD} X.
\]

\section{The set optimization perspective}
\label{SecSOP}

The constructions of the previous sections produce set-valued quantiles and a set-valued VaR in a natural way, but it might not be apparent how (much) these concepts are based on the complete lattice approach to  set optimization. It is the aim of this section to make this relationship transparent.
The basic reference is the survey \cite{HamelEtAl15Incoll}.

It is fundamental to introduce appropriate ``image spaces" for set-valued functions. In this note, lower quantiles and VaR map into
\[
\G(\R^d, C) = \cb{B \subseteq \R^d \mid B = \cl\co(B + C)}
\]
where $\cl$ denotes the topological closure, $\co$ the convex hull, and the addition $B + C = \cb{b+c \mid b \in B, c \in C}$ is the usual Minkowski addition of sets with the extension $B + \emptyset = \emptyset + B$ for all $B \in \G(\R^d, C)$. Thus, the addition in $\G(\R^d, C)$ has to be defined as $A \oplus B = \cl(A+B)$. The expression $A \ominus B$ is defined a parallel way. Together with a multiplication with non-negative reals defined by $s \cdot B = \cb{sb \mid b \in B}$ (in particular $s \cdot \emptyset = \emptyset$) for $s>0$ and $0 \cdot B = C$ (in particular $0 \cdot \emptyset = C$), the structure $(\G(\R^d, C), \oplus, \cdot)$ preserves as much of the structure of a linear space as possible (in \cite{HamelEtAl15Incoll} it is called a ``conlinear space"). Its order structure is even more important and summarized in the following result (see \cite{HamelEtAl15Incoll} and the references therein).

\begin{proposition}
\label{PropLattice}
The pair $(\G\of{\R^d, C}, \supseteq)$ is an order-complete lattice. If $\A \subseteq \G\of{\R^d, C}$, then
\[
\inf \A =  \cl\co\bigcup\limits_{A \in \A} A \quad \text{and} \quad
\sup \A = \bigcap\limits_{A \in \A} A
\]
where $\inf \A = \emptyset$ and $\sup \A = \R^d$ whenever $\A = \emptyset$. The greatest element in $(\G\of{\R^d, C}, \supseteq)$ is $\emptyset$, the least element is $\R^d$. 
\end{proposition}

``Order-complete" means that every subset has an infimum and a supremum. Remarkably, this is true without further assumptions to $C$ such as $\Int C \neq \emptyset$ or $(\R^d, \leq_C)$ is a vector lattice. Therefore, $C = \{0\}$ and $C = H^+(w)$ are valid options. The reader may observe that $(\G\of{\R^d, C}, \supseteq)$ shares its order features with $(\R\cup\cb{\pm\infty}, \leq)$ with the only exception that $\supseteq$ is not a total order.

Parallel, the set 
\[
\G(\R^d, -C) = \cb{B \subseteq \R^d \mid B = \cl\co(B - C)}
\]
is introduced with $B-C = \cb{b-c \mid b \in B, c \in C}$ and the same rules for $\emptyset$ as before as well as $B \ominus C = \cl(B - C)$. The pair $(\G\of{\R^d, -C}, \subseteq)$ is an order-complete lattice of ``downward" sets with the following formulas for infimum and supremum:
\[
\sup \A =  \cl\co\bigcup\limits_{A \in \A} A \quad \text{and} \quad
\inf \A = \bigcap\limits_{A \in \A} A.
\]
Note that the roles of union and intersection are swapped compared to the lattice $(\G\of{\R^d, C}, \supseteq)$ of ``upward" sets.

With these concepts in view, the lower $C$-quantile can be written as
\[
Q^-_{X,C}\of{p} = \cb{z \in \R^d \mid F_{X,C}(z) \geq p} = \inf\cb{z + C \mid z \in \R^d, \; F_{X,C}(z) \geq p}
\]
where the infimum on the right hand side now has to be taken in $(\G\of{\R^d, C}, \supseteq)$. The (closed convex) cone $C$ can be added by means of Proposition \ref{PropCQuantile} (a), hence the (very simple) function $z \mapsto z + C$ maps into $\G\of{\R^d, C}$. Seen in this way, the definition of the lower $C$-quantile is completely parallel to the definition of univariate lower quantiles: It is the $\G\of{\R^d, C}$-valued (lower) inverse of the function $z \to F_{X,C}(z)$.

The upper quantile function now becomes
\begin{align*}
Q^+_{X, w}\of{p} &  = \cb{z \in \R^d \mid P(w^TX < w^Tz) \leq p} \\
	& = \sup\cb{z - C \mid z \in \R^d, \; P(w^TX < w^Tz) \leq p}
\end{align*}
where the supremum is taken in $(\G\of{\R^d, -C}, \subseteq)$, thus $Q^+_{X, w}$ is the $\G\of{\R^d, -C}$-valued inverse of the function $z \mapsto P(w^TX < w^Tz)$.

Moreover, it might be observed that the definition of the $C$-distribution function and the lower $C$-quantile involve the scalar (!) infimum over the  family of distribution functions and the supremum in $(\G\of{\R^d, C}, \supseteq)$ (!) over the family of lower $w$-quantiles, thus, in this sense, they are also inverse to each other.

Consequently, $VaR_\alpha$ is a positively homogeneous, monotone and $\R^d$-translative $\G\of{\R^d, C}$-valued function (see Proposition \ref{PropVaR}). In contrast, the  VaRs defined by Embrechts/Puccetti in \cite{EmbrechtsPuccetti06JMA} as well as those in \cite{HamelHeyde10SIFIN} do not have convex values in general, hence they are much harder to handle when it comes to optimization, computation and in particular duality. For example, it is by no means clear how to define a multivariate AVaR starting from those definitions, but several options present themselves from the considerations above.

Finally, a remark concerning potential (risk) management applications might be in order. It has been claimed that a set-valued VaR such as the ones from \cite{EmbrechtsPuccetti06JMA} `can be unsuitable when we face real risk management problems' (\cite[p. 36]{CousinDiBernardino13JMA}). This point of view is shared by the authors of \cite[p. 112]{TorreLilloLanaido15IME} as they write `a multivariate VaR seen as a surface could bring problems with its interpretation.' We do not share this point of view. First, Definition \ref{DefVaR} produces a set-valued function which has a very clear (financial) interpretation. Secondly, in contrast to the mentioned references, we think that under a non-total order a ``single point risk measure" dupes a uniqueness property which is not inherent in the model: there always is an additional criterion according to which the single point is selected from a set (the Embrechts/Puccetti VaR in \cite{CousinDiBernardino13JMA}; a set obtained by shifting and turning the orthant $\R^d_+$ in a similar way as halfspaces are shifted and turned in Tukey's depth function in \cite{TorreLilloLanaido15IME}). When presented to a manager as ``the" risk compensating portfolio vector, (s)he might assume that this selection is the only choice; however,  there might be (and in general are) many more ``non-dominated" risk compensating portfolio vectors which might fit better if the manager has different weights for the components. As the formula in Proposition \ref{PropVaR} (a) shows, $VaR_\alpha(X)$ is robust with respect to the weights of the decision maker for the components of $z \in VaR_\alpha(X)$, but it provides flexibility for the management decision which is not present in the alternative approaches mentioned above. Simply put, one looses information if one selects only one point according to a fixed criterion instead of considering the whole set.

\section{Conclusions and perspectives}
\label{SecConPer}

We propose a `multidimensional counterpart of the quantiles of a random variable' (\cite[p. 1125]{BelloniWinkler11AS}) which are functions mapping into specific complete lattices of sets. This admits a calculus and applications parallel to the univariate case. In particular, quantile-based (financial) risk measures like the Value at Risk and stochastic orders can be introduced in a natural way. Our discussion also makes it desirable to investigate the following issues:

\begin{itemize}
\item to develop computational procedures for set-valued quantiles which can be based on ideas from computational geometry (see \cite{RousseeuwHubert15ArX} and the references herein) since for empirical distribution the method of choice would be the solution of linear vector optimization problems which in turn can be solved by tools closely related to computational convexity (see \cite{LoehneWeissing16EJOR} and the references therein),
\item to generalize the concepts to random cones $C$ with financial applications in view (the step from \cite{HamelHeyde10SIFIN} to \cite{HamelHeydeRudloff11MFE}),
\item to extend the approach to ``second order" constructions like the average or conditional value at risk and second order stochastic stochastic dominance,
\item to link the new concepts with dependence structures and study corresponding effects,
\item to study corresponding rank functions, outlyingness functions and similar concepts in the spirit of \cite{Serfling10JNS} and apply them to multivariate data analysis.
\end{itemize}

Finally, it might be a feasible attempt to deal with highly non-convex data sets via nonlinear ``scalarizations," i.e. one may replace the linear functions $z \mapsto w^Tz$ by particular classes of nonlinear ones. This has already been tried in \cite{HlubinkaEtAl10Kyb} in order to obtain ``weighted depth functions" which generalize Tukey's depth function, and this idea could be transferred to the context of this note in order to obtain the corresponding quantiles.

\bibliographystyle{plain}

\begin{thebibliography}{9999}

\bibitem{BelloniWinkler11AS} Belloni A, Winkler RL. On multivariate quantiles under partial orders. The Annals of Statistics 39(2):1125-79, 2011.

\bibitem{CascosMolchanov07FS} Cascos I, Molchanov I. Multivariate risks and depth-trimmed regions. Finance and Stochastics 11(3):373-97, 2007.

\bibitem{Chaudhuri96JASA} Chaudhuri, P. On a geometric notion of quantiles for multivariate data. Journal of the American Statistical Association 91(434):862-872.

\bibitem{CousinDiBernardino13JMA} Cousin A, Di Bernardino E. On multivariate extensions of Value-at-Risk. Journal Multivariate Analysis 119:32-46,  2013.

\bibitem{EmbrechtsPuccetti06JMA} Embrechts P, Puccetti G. Bounds for functions of multivariate risks. Journal Multivariate Analysis 97(2):526-47, 2006.

\bibitem{FoellmerSchied11Book3rd}
F\"{o}llmer H and Schied A. {\em Stochastic Finance: an Introduction in Discrete Time.} Walter de Gruyter Berlin New York, third revised and extended edition 2011.

\bibitem{HallinEtAl10AS} Hallin M, Paindaveine D, Siman M. Multivariate quantiles and multiple output regression quantiles: form $L_1$ optimization to halfspace depth. The Annals of Statistics 1:635-703, 2010

\bibitem{HamelHeyde10SIFIN} Hamel AH, Heyde F. Duality for set-valued measures of risk. SIAM Journal Financial Mathematics 1(1):66-95, 2010

\bibitem{HamelHeydeRudloff11MFE} Hamel AH, Heyde F, Rudloff, B. Set-valued risk measures for conical market models. 
Mathematics and Financial Economics, 5(1):1-28, 2011

\bibitem{HamelEtAl15Incoll} Hamel AH, Heyde F, L{\"o}hne A, Rudloff B, Schrage C. Set optimization--a rather short introduction.
In: Set optimization and applications--the state of the art. From set relations to set-valued risk measures. Springer Publishers Berlin 2015, pp. 65-141

\bibitem{HlubinkaEtAl10Kyb} Hlubinka, D, Kot{\'i}k, L, Venc{\'a}lek O. Weighted halfspace depth. Kybernetika 46(1):125-148, 2010

\bibitem{KongMizera12StSi} Kong L, Mizera I. Quantile tomography: using quantiles with multivariate data. Statistica Sinica 22(4):1589-1610, 2012

\bibitem{LoehneWeissing16EJOR} L\"ohne A, Wei{\ss}ing B. The vector linear program solver Bensolve -- notes on theoretical background.
European Journal Operational Research, http://dx.doi.org/10.1016/j.ejor.2016.02.039, 2016

\bibitem{MuellerStoyan02Book} M\"uller A, Stoyan D. {\em Comparison Methods for Stochastic Models and Risks.} John Wiley \& Sons, 2002.

\bibitem{RousseeuwHubert15ArX} Rousseeuw PJ, Hubert M. Statistical depth meets computational geometry: a short survey. arXiv preprint arXiv:1508.03828, 2015.

\bibitem{RousseeuwRuts99Met} Rousseeuw PJ, Ruts I. The depth function of a population distribution. Metrika 49(3):213-44,1999

\bibitem{SalvadoriMicheleDurante11HESS} Salvadori G, De Michele C, Durante F. On the return period and design in a multivariate framework. Hydrology and Earth Systems Sciences 15, 3293-3305, 2011

\bibitem{Serfling02SN} Serfling R. Quantile functions for multivariate analysis: approaches and applications. Statistica Neerlandica. 56(2):214-32, 2002.

\bibitem{Serfling06Incoll} Serfling, R. Depth functions in nonparametric multivariate inference, DIMACS Series in Discrete Mathematics and Theoretical Computer Science 72, pp. 1-16, 2006

\bibitem{Serfling10JNS} Serfling, R. Equivariance and invariance properties of multivariate quantile and related functions, and the role of standardisation. Journal of Nonparametric Statistics, 22(7):915-936, 2010

\bibitem{StruyfRousseeuv99JMA} Struyf AJ, Rousseeuw PJ. Halfspace depth and regression depth characterize the empirical distribution. Journal Multivariate Analysis 69(1): 135-53, 1999

\bibitem{TorreLilloLanaido15IME} Torres R, Lillo RE, Laniado H. A directional multivariate value at risk. Insurance: Mathematics and Economics. 65:111-23, 2015.

\bibitem{Tukey75} Tukey JW. Mathematics and the picturing of data. In: Proceedings of the International Congress of Mathematicians Vol. 2, pp. 523-531, 1975.

\bibitem{ZuoSerfling00AS} Zuo Y, Serfling R. General notions of statistical depth function. The Annals of Statistics 28(2):461-82, 2000.

\end{thebibliography}

\end{document}